\documentclass[12pt]{article}
\usepackage{amsmath,amsxtra,amssymb,latexsym, amscd,amsthm}

\advance\voffset-1truecm\relax
\advance\hoffset-1truecm\relax

\def\R{{\Bbb R}}
\def\N{{\Bbb N}}
\def\R{{\Bbb R}}
\def\C{{\Bbb C}}
\def\P{{\Bbb P}}

\def\Z{{\Bbb Z}}

\numberwithin{equation}{section}

\newcommand {\Cal}{\mathcal}

\newtheorem{lemma}{Lemma}[section]

\newtheorem{corollary}[lemma]{Corollary}

\newtheorem{proposition}[lemma]{Proposition}
\newtheorem{definition}[lemma]{Definition}
\newtheorem{remark}[lemma]{Remark}
\def\leq{\leqslant}

\begin{document}
\title{A SECOND MAIN THEOREM FOR\\
MOVING HYPERSURFACE TARGETS}

\author{$\quad$Gerd Dethloff and Tran Van Tan}
\date{$\quad$}
\maketitle
\vspace{-0.5cm}
\begin{abstract}
\noindent In this paper, we prove a Second Main Theorem for  algebraically
nondegenerate meromorphic maps of $\C^m$ into $\C P^n$ and slowly
moving hypersurfaces targets $Q_j \subset\C P^n,$ $j=1,\dots,q
\;(q\geq n+2)$  in (weakly) general position. This generalizes the Second Main Theorem for fixed hypersurface targets in general position, obtained by M. Ru in \cite{b13}. We also introduce a truncation, with an explicit estimate of the truncation level, into this Second
Main Theorem with moving targets, thus generalizing the main result 
of An-Phuong \cite{AP}.
\end{abstract}

\section{Introduction}
For $z = (z_1,\dots,z_m) \in \C^m$, we set
$\Vert z \Vert = \Big(\sum\limits_{j=1}^m |z_j|^2\Big)^{1/2}$
and define
\begin{align*}
B(r) &= \{ z \in \C^m : \Vert z \Vert < r\},\quad
S(r) = \{ z \in \C^m : \Vert z \Vert = r\},\\
d^c &= \dfrac{\sqrt{-1}}{4\pi}(\overline \partial - \partial),\quad
\Cal V = \big(dd^c \Vert z\Vert^2\big)^{m-1},\; \sigma = d^c
\text{log}\Vert z\Vert^2 \land \big(dd^c\text{log}\Vert z
\Vert\big)^{m-1}.
\end{align*}

\hbox to 5cm {\hrulefill }

 \noindent {\small Mathematics Subject
Classification 2000: Primary 32H30; Secondary 32H04, 32H25,
14J70.}

\noindent {\small Key words and phrases:  Nevanlinna theory,
Second Main Theorem.}

 {\small The first named author was partially supported by the Fields Institute Toronto. The second named author was partially supported by the
post-doctoral research program of the Abdus Salam International Centre for Theoretical Physics.}

Let $L$ be a positive integer or $+\infty$ and $\nu$ be a divisor on
$\C^m.$ Set $ |\nu| = \overline {\{z : \nu(z) \neq 0\}}.$
 We define the counting function of $\nu$ by

\begin{align*}
N^{(L)}_\nu(r) := \int\limits_1^r \frac{n^{(L)}(t)}{t^{2m-1}}dt\quad
(1 < r < +\infty),
\end{align*}
where
\begin{align*}
n^{(L)}(t) &= \int\limits_{|\nu | \cap B(t)} \text{min}\{\nu
,L\}\cdot \Cal V\ \quad
\text{for}\quad m \geq 2 \ \text{and}\\
n^{(L)}(t) &= \sum_{|z| \leq t}\text{min}\{ \nu(z),L\} \qquad\quad
\text{for}\quad m = 1.
\end{align*}

Let $F$ be a nonzero holomorphic function on $\C^m$. For a set
$\alpha = (\alpha_1,\dots,\alpha_m)$ of nonnegative integers, we set
$|\alpha| := \alpha_1 + \dots + \alpha_m$ and $D^\alpha F :=
\dfrac{\partial^{|\alpha|}} {\partial^{\alpha_1}z_1 \cdots
\partial^{\alpha_m}z_m}\,\cdotp$
 We define the zero divisor  $\nu_F$ of $F$ by
\begin{align*}
\nu_F(z) = \max \big\{ p : D^\alpha F(z) = 0 \ \text{for all
$\alpha$ with}\ |\alpha| < p \big\}.
\end{align*}

Let $\varphi$ be a nonzero meromorphic function on $\C^m$. The zero
divisor $\nu_\varphi$ of $\varphi$ is defined as follows: For each
$a \in \C^m$, we choose nonzero holomorphic functions $F$ and $G$ on
a neighborhood $U$ of $a$ such that $\varphi = \dfrac{F}{G}$ on $U$
and $\text{dim}\big(F^{-1}(0) \cap G^{-1}(0)\big) \leq m-2$, then we
put $\nu_\varphi(a) := \nu_F(a)$.

Set $N_\varphi^{(L)}(r):=N_{\nu_\varphi}^{(L)}(r).$ For brevity we
will omit the character ${}^{(L)}$ in the counting function if
$L=+\infty.$

Let $f$ be a meromorphic map of $\C^m$ into $\C P^n$. For arbitrary
fixed homogeneous coordinates $(w_0: \cdots : w_n)$ of $\C P^n$, we
take a reduced representation $f = (f_0 : \cdots : f_n)$, which
means that each $f_i$ is a holomorphic function on $\C^m$ and $f(z)
= (f_0(z) : \cdots : f_n(z))$ outside the analytic set $\{ z :
f_0(z) = \cdots = f_n(z) = 0\}$ of codimension $\geq 2$. Set $\Vert
f  \Vert = \max \{ |f_0|, \dots  , |f_n| \}$.

The characteristic function of $f$ is defined by
\begin{align*}
T_f(r) := \int\limits_{S(r)}\text{log}\Vert f \Vert \sigma -
\int\limits_{S(1)} \text{log}\Vert f \Vert \sigma ,\quad 1 < r < +\infty.
\end{align*}

 For a meromorphic function $\varphi$ on $\C^m$, the characteristic function
$T_\varphi(r)$ of $\varphi$ is defined by considering  $\varphi$ as
a meromorphic map of $\C^m$ into $\C P^1$. We have the following
Jensen's formula :
\begin{equation*}
N_{\varphi }(r)-N_{\frac{1}{\varphi }}(r)=\int\limits_{S(r)}\text{log}%
|\varphi |\sigma -\int\limits_{S(1)}\text{log}|\varphi |\sigma .
\end{equation*}

Let $f$ be a nonconstant meromorphic map of $\C^m$ into $\C P^n$.
We say that a meromorphic function $\varphi$ on $\C^m$ is ``small"
with respect to $f$ if $T_\varphi(r) = o(T_f(r))$ as $r \to \infty$
(outside a set of finite Lebesgue measure).

Denote by $\Cal K_f$ the set of all ``small" (with respect to $f$)
meromorphic functions on $\C^m$. By Theorem 5.2.29 of \cite{b20} or by Corollary 5.7 in \cite{b7} we
easily get that any rational expression of functions in $\Cal K_f$ is
still ``small'' (with respect to $f$), in particular $\Cal K_f$ is a field.

For a homogeneous polynomial $Q \in \Cal K_f[x_0,\dots,x_n]$ of
degree $d \geq 1$ with $Q(f_0, \dots , f_n) \not\equiv 0$, we define
$$
N^{(L)}_f(r,Q) :=
N^{(L)}_{Q(f_0,\dots,f_n)}(r)\:\text{and}\;\delta_f (Q) = \lim_{r
\to \infty}\inf \Big(1 - \frac{N_f(r,Q)}{d \cdot T_f(r)}\Big).$$

Denote by $Q(z)$ the homogeneous  polynomial over $\C$ obtained by
evaluating the coefficients of $Q$ at a specific point $z \in \C^m$
in which all coefficient functions of $Q$ are holomorphic.

For a positive integer $d$, we set
\begin{align*}
\Cal T_d := \big\{ (i_0,\dots,i_n) \in \N_0^{n+1} :
i_0 + \dots + i_n = d \big\}.
\end{align*}
Let
\begin{align*}
Q_j = \sum\limits_{I \in \Cal T_{d_j}} a_{jI}x^I \quad (j = 1,\dots,q)
\end{align*}
 be homogeneous polynomials in $\Cal K_f[x_0,\dots,x_n]$ with
$\text{deg}\,Q_j = d_j \geq 1$, where $x^I = x_0^{i_0} \cdots x_n^{i_n}$
for $x = (x_0,\dots,x_n)$ and $I = (i_0, \dots,i_n)$.
Denote by $\Cal K_{\{Q_j\}_{j=1}^q}$ the field over $\C$ of   all
meromorphic functions on $\C^m$ generated by
$\big\{ a_{jI} : I \in \Cal T_{d_j}, j \in \{1,\dots,q\}\big\}$.
It is clearly a subfield of $\Cal K_f$.
Denote by
$\tilde{\Cal K}_{\{Q_j\}_{j=1}^q} \subset \Cal K_{\{Q_j\}_{j=1}^q}$
the subfield generated by all quotients
$\big\{\frac{a_{jI_{1}}}{a_{jI_{2}}} :a_{jI_{2}}\not=0,  I_{1}, I_{2}
 \in \Cal T_{d_j}; j \in \{1,\dots,q\} \big\}$.
We say that $f$ is algebraically nondegenerate over
$\Cal K_{\{Q_j\}_{j=1}^q}$
(respectively over $\tilde{\Cal K}_{\{Q_j\}_{j=1}^q}$)
if there is no nonzero homogeneous polynomial
$Q \in \Cal K_{\{Q_j\}_{j=1}^q}[x_0,\dots,x_n]$
(respectively  $Q \in \tilde{\Cal K}_{\{Q_j\}_{j=1}^q}[x_0,\dots,x_n]$)
such that
$Q(f_0,\dots,$ $f_n) \equiv 0$.

We say that a set $\{Q_j\}_{j=1}^q$ $(q \geq n+1)$ of homogeneous
polynomials in $\Cal K_f [x_0,\dots,$ $x_n]$ is admissible (or in (weakly) general position)
if there exists $z \in \C^m$
in which all coefficient functions of all $Q_j$, $j=1,...,q$ are holomorphic and such that for any
$1 \leq j_0 < \dots < j_n \leq q$ the system of equations
\begin{align} \label{zz}
\left\{ \begin{matrix}
Q_{j_i}(z)(x_0,\dots,x_n) = 0\cr
0 \leq i \leq n\end{matrix}\right.
\end{align}
has only the trivial solution $(x_0, \dots , x_n) = (0,\dots,0)$ in $\C^{n+1}$.
We remark that in this case this is true for the generic $z \in \C^m$.

As usual, by the notation ``$\Vert P$" we mean the assertion $P$ holds
for all $r \in [1, +\infty)$ excluding a Borel subset $E$ of $(1, +\infty)$
with $\displaystyle{\int\limits_E} dr < +  \infty$.

\noindent{\bf Main Theorem}. {\it Let $f$ be a nonconstant
meromorphic map of $\C^m$ into $\C P^n$. Let
$\big\{Q_j\big\}_{j=1}^q$ be an admissible set of homogeneous
polynomials in $\Cal K_f [x_0,\dots,x_n]$ with $\deg Q_j = d_j \geq
1$. Assume that $f$ is algebraically nondegenerate over $\tilde
{\Cal K}_{\{Q_j\}_{j=1}^q}$. Then for any $\varepsilon > 0,$ there
exist  positive integers $L_j\;(j=1,\dots,q)$, depending only on
$n,$ $\varepsilon$ and $d_j \;(j=1,\dots, q)$ in an explicit way
such that}
\begin{align*}
\Vert (q-n-1-\varepsilon) T_f(r) \leq \sum_{j=1}^q \frac{1}{d_j}
N^{(L_j)}_f(r,Q_j).
\end{align*}

We note that, for fixed hypersurface targets, in 1979, Shiffman
\cite{b16} conjectured that if $f$ is an algebraically nondegenerate
holomorphic map of $\C$ into $\C P^n$ and $D_1,\cdots,D_q$ are
hypersurfaces in $\C P^n$ in general position, then $\sum_{j=1}^q
\delta_f (D_j) \leq n+1.$ This conjecture was proved by  Ru
\cite{b13} in 2004, and recently even generalized by Ru \cite{Ru} to
fixed hypersurface sections of projective varieties in general position.  As a corollary of the Main Theorem we get the
generalization of his result in \cite{b13} for moving
targets.
\begin{corollary} [Shiffman conjecture for moving
hypersurfaces]  Under the same assumption  as in the Main theorem,
we have
\begin{align*}
\sum_{j=1}^q \delta_f (Q_j) \leq n+1.
\end{align*}
\end{corollary}

We also note that for the case of moving hyperplanes $(d_1 = \dots =
d_q = 1)$, and multiplicities which are not truncated,
 the above theorem was first proved by  Ru and  Stoll in 1991 \cite{b14}. In 2002,  Tu \cite{tu} introduced a truncation into the Second
 Main Theorem of Ru-Stoll, but the truncation level is is not estimated.
Furthermore, after the first version \cite{DT} of this paper was published, in which the truncation level was not estimated, neither, An-Phuong \cite{AP} gave a truncation with an explicit
estimate for the Second Main Theorem for fixed hypersurfaces. So
our Main Theorem, now also with an explicit estimate of the truncation level, is also a generalization of their result to moving
hypersurfaces. In the special case of fixed hypersurfaces our
estimate for the truncation is still slightly better, but,
at least in the case when all hypersurfaces are of the same degree, still of the same order than theirs.

\begin{proposition} With the notation of  our Main Theorem, we have   $$L_j\leq
\frac{d_j\cdot\binom{n+N}{n}t_{p_0+1}-d_j}{d}+1\:,$$ where $d$ is the
least common multiple of the $d_j$'s and
$$N=d\cdot[2(n+1)(2^n-1)(nd+1)\epsilon^{-1}+n+1]\:,$$
$$p_0=\big[\frac{\big( \binom{n+N}{n}^2.\binom{q}{n}-1
\big).\log\big(
\binom{n+N}{n}^2.\binom{q}{n}\big)}{\log(1+\frac{\epsilon}{2\binom{n+N}{n}N})}+1\big]^2,$$
$$\text{and}\quad
 t_{p_0+1}
<\Bigg(\binom{n+N}{n}^2.\binom{q}{n}+p_0\Bigg)^{
\binom{n+N}{n}^2.\binom{q}{n}-1},$$ where we denote $[x]:=\max\{k\in
\Z: k\leq x\}$ for a real number $x.$ Furthermore, in the case of
fixed hypersurfaces $(Q_j\in\C[x_0,\dots,x_n], j=1,\dots, q),$ we
have $t_p=1$ for all positive integers $p$, so we get a better estimate:
$$L_j\leq\frac{d_j\cdot\binom{n+N}{n}-d_j}{d}+1.$$
\end{proposition}

 \begin{remark} The Main Theorem holds, more generally,
for an admissible set of polynomials $\big\{Q_j\big\}_{j=1}^q$ such
that only the quotients $\big\{\frac{a_{jI_{1}}}{a_{jI_{2}}}
:a_{jI_{2}}\not=0,  I_{1}, I_{2}
 \in \Cal T_{d_j}; j \in \{1,\dots,q\} \big\}$ lie in $\Cal K_f$, under the
 condition that one replaces the $N^{(L)}_f(r,Q_j)$ by
$N^{(L)}_f(r,\frac{1}{a_{jI_{2}}}Q_j)$, where $a_{jI_{2}}\not=0$ can
be any nonzero coefficient of $Q_{j}$, $j=1,...,q$. This follows
immediately from the Main Theorem, applied to the set of polynomials
$\big\{\tilde Q_j \big\}_{j=1}^q$, where $\tilde Q_j
:=\frac{1}{a_{jI_{2}}}Q_j $. For more details, see the beginning of
section 4.
\end{remark}

The proof of our Main Theorem (including the one of Proposition 1.2) consists of three main parts, in which the second and the third one
are considerably  more complicated than this can be done for fixed
hypersurfaces with their notion of general position:

In the first part (chapter 4 until equation (\ref{eq37})) we use the idea of Corvaja-Zannier \cite{Co} and Ru \cite{b13}
to estimate $log\prod_{j=1}^q\vert Q_j(f)\vert.$ However, we have
to pass many difficulties which come both from the facts that the concept ``in general position'' in our paper is more general than in
Corvaja-Zannier's and Ru's paper  and that the field
$\Cal K_f$ is not algebraically closed in general, so we cannot use
any more Hilbert's Nullstellensatz.
 Instead we have to
use explicit results on resultants respectively discriminant
varieties for universal families of configurations of $q$
hypersurfaces in $\C P^n$, generalizing, among others, 
considerably Hilbert's Nullstellensatz (see \cite{b10}, chapter IX).
This allows us to deal with such hypersurfaces with ``variable'' coefficients, namely in $\Cal K_f$,
but by specialization to the fibers to have nevertheless complex solutions of these configurations of hypersurfaces. Another problem
related to the fact that $\Cal K_f$ is not algebraically closed in general
is that the proof of the fact that admissible families of polynomials
in $\Cal K_f[x_0,...,x_n]$ give regular families does not follow
any more directly from Hilbert's Nullstellensatz, but needs another
time resultants, as well as results on parameter systems in Cohen-Macauley rings.

 In the second part (up to equation (\ref{13}), we estimate the ``error term'' of equation (\ref{eq37}), relating it moreover to a Wronskian, which will
 become crucial to give the truncation in the third part. It is in
 particular here where
 generalizing the coefficients from constants to
meromorphic functions (although slowly growing ones) complicates
substantially the analysis, especially with respect to the
Wronskians and the Lemma of Logarithmic Derivative. Here we have to
introduce technics known from Value Distribution Theory of moving
hyperplanes (which we take from Shirosaki \cite{b17}), and to adopt
them from the hyperplane to the hypersurface case. Another
complication compared to the moving hyperplane case is that we
cannot use once and for all reduced representations for the
coefficient functions of the polynomials giving the moving
hypersurfaces, which needs a special care while we take pointwise
maxima or minima of their norms and while we estimate error terms.
It is only at the end of the proof when we use a Lemma of
Logarithmic Derivative for wronskians,  where we pass to
a reduced representation of a particular meromorphic map from $\C^m$
with monomial coefficients in the components of $f$ and the
coefficients of the $Q_j, j=1, \dots q$. We finally remark that in this part,
instead using the lemma of logarithmic derivative, we also could have
used Theorem 2.3 of Ru \cite{b12}.

In the third part, truncation is obtained. Here the
concept "resultants of  homogenous polynomials" and Wronskians are
used again, now to estimate the corresponding divisors. The use of this tool, which is not necessary in the case of fixed hypersurfaces, is
necessary 
in the case of moving hypersurfaces because of our very general notion of general position, in order to control what happens over
the divisor where the resultant vanishes, this means where the
hypersurfaces are not in general position.  

We finally remark that we prefered to prove our result right away for
meromorphic maps rather than only for the most important special case, namely entire holomorphic curves, 
since this proof is only around two pages longer than the one we could have given for entire curves.

\section{Some lemmas}
We first recall some classical results on resultants,
see Lang \cite{b10}, section IX.3, for the precise definition, the existence and for the principal properties of resultants, as well as Eremenko-Sodin \cite{b4}, page 127:
Let $\big\{Q_j\big\}_{j=0}^n$ be a set of homogeneous
polynomials  of common degree $d \geq 1$ in
$\Cal K_f[x_0,\dots,x_n]$
 \begin{align*}
Q_j = \sum_{I \in \Cal T_d} a_{jI}x^I,\quad a_{jI} \in \Cal K_f \quad
(j = 0,\dots,n).
\end{align*}
Let $T = (\dots,t_{kI},\dots)$ \ ($k \in \{0,\dots,n\}$, $I \in \Cal T_d$)
be a family of variables. Set
\begin{align*}
\widetilde Q_j = \sum_{I \in \Cal T_d} t_{jI}x^I \in \Z[T,x],\quad
j = 0,\dots, n.
\end{align*}
Let $\widetilde R \in \Z[T]$ be the resultant of $\widetilde Q_0, \dots,
\widetilde Q_n$. This is a polynomial in the variables 
$T = (\dots,t_{kI},\dots)$ \ ($k \in \{0,\dots,n\}$, $I \in \Cal T_d$)
with integer coefficients, such that the condition 
$\widetilde R (T) =0$ is necessary and sufficient for the
existence of a nontrivial solution 
$(x_0, \dots , x_n) \not= (0,\dots,0)$ in $\C^{n+1}$
of the system of equations
\begin{align} \label{z}
\left\{ \begin{matrix}
\widetilde Q_{j}(T)(x_0,\dots,x_n) = 0\cr
0 \leq i \leq n\end{matrix}\right. \:.
\end{align}
From equations (\ref{z}) and (\ref{zz}) is follows immediately
that if $$\big\{Q_j= \widetilde Q_j(a_{jI})(x_0, \dots, x_n)\,, \:j=0, \dots , n\big\}$$ is an admissible set,
\begin{equation}R := \widetilde R(\dots, a_{kI}, \dots) \not\equiv 0\,.\label{zzz}
\end{equation}
Furthermore, since  $a_{kI} \in \Cal K_f$, we have  
$R \in \Cal K_f$.
 We finally will need the following result on resultants,
 which is contained in Theorem 3.4 in \cite{b10} (see also Eremenko-Sodin \cite{b4}, page 127, for a similar result):
\begin{proposition}\label{lang} There exists a positive integer $s$
and polynomials $\big\{\widetilde
b_{ij}\big\}_{0 \leq i, j \leq n}$ in $\Z[T,x]$, which are (without loss of generality) zero or
homogenous in $x$ of degree $s-d$,
such that
\begin{align*}
x_i^s \cdot \widetilde R = \sum_{j=0}^n \widetilde b_{ij} \widetilde Q_j\quad
\text{for all}\ i \in \{0,\dots,n\}.
\end{align*}
  \end{proposition}

Let $f$ be a nonconstant meromorphic map of $\C^m$ into $\C P^n$.
Denote by $\Cal C_f$ the set of
all    non-negative functions
$h : \C^m \setminus A \longrightarrow [0,+\infty] \subset \overline{\R}$,
which are of the form
\begin{equation}\label{expr}
\frac{|g_1|+ \dots + |g_k|}{|g_{k+1}|+ \dots + |g_l|}\: ,
\end{equation}
where $k,l \in \N$, $g_1, \dots , g_l \in \Cal K_f \setminus \{0\}$ and
$A \subset \C^m$, which may depend on $g_1, \cdots , g_l\,$,
is an analytic set of codimension at least two.
By Jensen's formula and the First Main Theorem we have
\begin{align*}
\int\limits_{S(r)} \text{log} |\phi | \sigma = o(T_f(r))\quad
\text{as}\ r \to \infty
\end{align*}
for $\phi \in \Cal K_f \setminus \{0\}$.
Hence,  for any $h \in \Cal C_f$, we have
\begin{align*}
\int\limits_{S(r)} \text{log} h \sigma = o(T_f(r))\quad
\text{as}\ r \to \infty.
\end{align*}
It is easy to see that sums, products and quotients of functions
in $\Cal C_f$ are again in $\Cal C_f$.  We would like to
point out that, in return, given any functions
$g_1, \cdots , g_l \in \Cal K_f \setminus \{ 0 \}$, any expression of the form (\ref{expr})
is in fact a well defined function (with values in $[0,+\infty]$) outside  an analytic subset $A$ of
codimension at least two, even though all the
$g_1, \cdots , g_l $ can have common pole or zero divisors in codimension one.
\begin{lemma}\label{lem1}
Let $\big\{Q_j\big\}_{j=0}^n$ be a set of homogeneous
polynomials  of degree $d$ in $\Cal K_f[x_0,\dots,x_n]$. Then there exists a
function $h_1 \in \Cal C_f$ such that, outside an analytic set of
$\C^m$ of codimension at least two,
\begin{align*}
\max_{j \in \{0,\dots,n\}}
|Q_j(f_0,\dots,f_n)| \leq h_1 \cdot \Vert f \Vert^d.
\end{align*}
If, moreover, this set of homogeneous polynomials is admissible, then
there exists a nonzero
function $h_2 \in \Cal C_f$ such that, outside an analytic set of
$\C^m$ of codimension at least two,
\begin{align*}
h_2 \cdot \Vert f \Vert^d \leq
\max_{j \in \{0,\dots,n\}}
|Q_j(f_0,\dots,f_n)| .
\end{align*}

\end{lemma}

\begin{proof}
Assume that  \begin{align*}
Q_j = \sum_{I \in \Cal T_d} a_{jI}x^I,\quad a_{jI} \in \Cal K_f \quad
(j = 0,\dots,n).
\end{align*}
We have, outside a proper analytic set of $\C^m$,
\begin{align}
|Q_j(f_0,\dots,f_n)| = \Big|\sum_{I \in \Cal T_d}a_{jI}f^I\Big|
\leq \sum_{I \in \Cal T_d}|a_{jI}|\cdot \Vert f\Vert^d. \label{eq21}
\end{align}
Set
\begin{align*}
h_1 := \sum_{j=0}^n\sum_{I \in \Cal T_d}|a_{jI}|.
\end{align*}
Then $h_1 \in \Cal C_f$, since $a_{jI} \in \Cal K_f$ \
($j \in \{0, \dots,n\}$,  $I \in \Cal T_d$).
By (\ref{eq21}), we get
\begin{align*}
|Q_j (f_0,\dots,f_n)| \leq |h_1| \cdot \Vert f\Vert^d \quad
\text{for all}\ j \in \{0,\dots,n\}.
\end{align*}
So we have
\begin{align}
\max_{j \in \{0,\dots,n\}}|Q_j(f_0,\dots,f_n)| \leq h_1 \cdot
\Vert f \Vert^d . \label{eq22}
\end{align}
All expressions in the last inequality are well defined and continuous
(as functions with values in $[0,+ \infty]$) outside  analytic sets
of codimension at least two.
Since $\Vert f \Vert^d$ is a real-valued function which  is zero
only on an analytic subset of $\C^m$ of codimension at least two,
 this inequality still holds outside an analytic subset of $\C^m$ of codimension at least two.

In order to prove the second inequality, by Proposition \ref{lang}
and its notations we have:
There exists a positive integer $s$
and polynomials $\big\{\widetilde
b_{ij}\big\}_{0 \leq i, j \leq n}$ in $\Z[T,x]$,
zero or
homogenous in $x$ of degree $s-d$,
such that
\begin{align*}
x_i^s \cdot \widetilde R = \sum_{j=0}^n \widetilde b_{ij} \widetilde Q_j\quad
\text{for all}\ i \in \{0,\dots,n\}.
\end{align*}
Moreover,
$R = \widetilde R(\dots, a_{kI}, \dots) \not\equiv 0$.
Set
\begin{align*}
b_{ij} = \widetilde b_{ij}\big((\dots,a_{kI},\dots), (f_0,\dots,f_n)\big),\quad
0 \leq i, j \leq n.
\end{align*}
Then, we get

\begin{align*}
f_i^s \cdot R = \sum_{j=0}^n b_{ij} \cdot Q_j(f_0,\dots,f_n)\quad
\text{for all}\ i \in \{0,\dots,n\}.
\end{align*}
So we have, outside a proper analytic set of $\C^m$:
\begin{align}
|f_i^s \cdot R| &= \Big|\sum_{j=0}^n b_{ij} \cdot
Q_j(f_0,\dots,f_n)\Big| \notag\\
&\leq \sum_{j=0}^n |b_{ij}| \cdot \max_{k \in \{0,\dots,n\}}
|Q_k(f_0,\dots,f_n)| \label{eq23}
\end{align}
for all $i \in \{0,\dots,n\}$.
We write
\begin{align*}
b_{ij} = \sum_{I \in \Cal T_{s-d}} \gamma_I^{ij} f^I,\quad
\gamma_I^{ij} \in \Cal K_f.
\end{align*}
By (\ref{eq23}), we get
\begin{align*}
|f_i^s \cdot R| \leq \sum_{\begin{matrix} \scriptstyle{0 \leq j \leq n}\cr
\noalign{\vskip-0.15cm}
\scriptstyle{I \in \Cal T_{s-d}}\end{matrix}}
\big|\gamma_I^{ij}\big| \cdot \Vert f \Vert^{s-d} \cdot
\max_{k \in \{0,\dots,n\}} |Q_k(f_0,\dots,f_n)|,\quad
i \in \{ 0,\dots,n\}.
\end{align*}
So
\begin{align}
\frac{|f_i|^s}{\Vert f \Vert^{s-d}} \leq \sum_{\begin{matrix}
\scriptstyle{0 \leq j \leq n}\cr
\noalign{\vskip-0.15cm}
\scriptstyle{I \in \Cal T_{s-d}}\end{matrix}} \Big|\frac{\gamma_I^{ij}}
{R}\Big| \max_{k \in \{0,\dots,n\}} |Q_k(f_0,\dots,f_n)| \label{eq24}
\end{align}
for all $i \in \{0,\dots,n\}$.
Set
\begin{align*}
h_2 = \frac{1}{\sum_{i=0}^n \sum_{\begin{matrix} \scriptstyle{0 \leq j \leq n}\cr
\noalign{\vskip-0.15cm}
\scriptstyle{I \in \Cal T_{s-d}}\end{matrix}}
\Big|\frac{\gamma_I^{ij}}{R}\Big|}.
\end{align*}
Then $h_2 \in {\Cal C}_f$, since $\gamma_I^{ij}, R \in {\Cal K}_f$ and
$R \not\equiv 0$.
By (\ref{eq24}) and since $\Vert f\Vert $ was the maximum norm,
so $\Vert f\Vert = |f_i|$ for some $i=0, \dots , n$ (which may depend on $z \in \C^m$), we have
\begin{align}
h_2 \cdot \Vert f \Vert^d \leq  \max_{j \in \{0,\dots,n\}}
|Q_j(f_0,\dots,f_n)|\:. \label{eq25}
\end{align}
By (\ref{eq22}) and (\ref{eq25}) and by the same observations as for
the first inequality we get Lemma \ref{lem1}
\end{proof}

Consider  meromorphic functions  $F_0,\dots,F_n$ on $\C^m$,
and put $F = (F_0,\dots,F_n)$.
For each $a \in \C^m$, we denote by $\Cal M _a$ the field of all germs of
meromorphic functions on $\C^m$ at $a$ and, for $p=1,2,\dots$ by $\Cal F^p$
the $\Cal M_a$-sub vector space of $\Cal M_a^{n+1}$ which is generated by
the set $\{D^\alpha F := (D^\alpha F_0,\dots,D^\alpha F_n) :
|\alpha|\leq p\}$. Set $\ell_F(p) = \text{dim}_{\Cal M_a}\Cal F^p$,
which does not depend on  $a \in \C^m$. As a general reference for this construction
and for the following definition, see \cite{b5} and \cite{b6}.

\begin{definition}\label{lem4} {\rm (see \cite{b6}, Definition 2.10)
Assume  that meromorphic functions $F_0,\dots,F_n$ on $\C^m$ are linearly
independent over $\C$. For $(n+1)$ vectors
$\alpha^i = (\alpha_{i1},\dots,\alpha_{im})$ $(0 \leq i \leq n)$ composed
of nonnegative integers $\alpha_{ij}$, we call a set $\alpha = (\alpha^0,\dots,
\alpha^n)$ an admissible  set for $F := (F_0,\dots,F_n)$ if
$$\big\{ D^{\alpha^0}F, \dots,D^{\alpha^{\ell_{F}(p)-1}}F \big\}$$
is a basis of $\Cal F^p$ for each $p = 1, 2, \dots$,
$p_0 := \min \{ p' : \ell_{F}(p') = n+1\}$.}
\end{definition}

By definition, for an admissible set $\alpha = (\alpha^0,\dots,\alpha^n)$ for
$F = (F_0,\dots,F_n)$ we have
\begin{align*}
W^\alpha (F_0,\dots,F_n) :=
\text{det}\,(D^{\alpha^0}F, \dots,D^{\alpha^n}F) \not\equiv 0.
\end{align*}

\begin{lemma}\label{lem5} {\rm (\cite{b6}, Proposition 2.11)}
For arbitrarily given linearly independent meromorphic functions
$F_0, \dots, F_n$ on $\C^m$, there exists an admissible set
$\alpha = (\alpha^0,\dots,\alpha^n)$ with
 $|\alpha|:= \sum_{i=0}^{n}|\alpha^{i}| \leq
\dfrac{n(n+1)}{2}$.
\end{lemma}

\begin{lemma}\label{lemz}
For arbitrarily given linearly independent meromorphic functions
$F_0, \dots, F_n$ on $\C^m$,
$$p_0 := \min \{ p' : \ell_{F}(p') = n+1\} \leq n \,.$$
\end{lemma}

\begin{proof}
This is an easy corollary of Fujimoto \cite{b6}, Proposition 2.9,
since $F$ is at least of rank one, or of Fujimoto \cite{b5}, Proposition 4.5.
\end{proof}

\begin{lemma}\label{lem6} {\rm (generalization of \cite{b6}, Proposition 2.12)}
Let $\alpha = (\alpha^0, \dots,\alpha^n)$ be an admissible set for
$F=(F_0,\dots,F_n)$ and let $h$ be a nonzero meromorphic
function on $\C^m$.
Then
\begin{align*}
W^\alpha (hF_0,\dots,hF_n) = h^{n+1}W^\alpha(F_0,\dots,F_n).
\end{align*}
\end{lemma}

\begin{proof} For holomorphic functions $h$ this is Proposition 2.11 in \cite{b6},
and its proof argument still holds for holomorphic functions defined only
on a Zariski open subset of $\C^m$. Hence, the case of a meromorphic  $h$
follows by the identity theorem.
\end{proof}

We also will need the following variant of the logarithmic derivative lemma:

\begin{lemma}\label{lem7}
Let $f$ be a linearly nondegenerate meromorphic map of $\C^m$ into
$\C P^n$ with reduced representation $f = (f_0,\dots,f_n)$.
Let $\alpha = (\alpha^0, \dots,\alpha^n)$ be an admissible
 set for $(f_0,\dots,f_n)$.  Then
\begin{align*}
\Big\Vert \int\limits_{S(r)} \text{\rm log}^+
\Big|\frac{W^\alpha (f_0,\dots,f_n)}{f_0 \cdots f_n}\Big| \sigma =
o(T_f(r)).
\end{align*}
\end{lemma}

\begin{proof}
By Lemma \ref{lem6} we have
\begin{align*}
\int\limits_{S(r)} \text{log}^+&\Big|\frac{W^\alpha (f_0,\dots,f_n)}
{f_0 \cdots f_n}\Big| \sigma = \int\limits_{S(r)} \text{log}^+
\Big|\frac{W^\alpha \Big(1, \dfrac{f_1}{f_0},\dots, \dfrac{f_n}{f_0}\Big)}
{1 \cdot \dfrac{f_1}{f_0} \cdots \dfrac{f_n}{f_0}}\Big| \sigma \\
&\leq \int\limits_{S(r)} \Big(K_1 \sum_{\begin{matrix}
\scriptstyle{0 \leq i \leq n}\cr
\noalign{\vskip-0.15cm}
\scriptstyle{1 \leq j \leq n}\end{matrix}} \text{log}^+
\Bigg|\frac{D^{\alpha^i}\Big(\dfrac{f_j}{f_0}\Big)}
{\dfrac{f_j}{f_0}}\Bigg| + K_2 \Big)\sigma\\
&\leq K_1 \sum_{\begin{matrix}
\scriptstyle{0 \leq i \leq n}\cr
\noalign{\vskip-0.15cm}
\scriptstyle{1 \leq j \leq n}\end{matrix}} \
\int\limits_{S(r)} \text{log}^+
\Bigg| \frac{D^{\alpha^i}\Big(\dfrac{f_j}{f_0}\Big)}
{\dfrac{f_j}{f_0}}\Bigg|\sigma + K_3\:,
\end{align*}
where $K_1$, $K_2$, $K_3$ are constant not depending on $r$.
On the other hand, by Theorem 2.6 in \cite{b6}, we have
\begin{align*}
\Big\Vert \int\limits_{S(r)} \text{log}^+ \Bigg|\frac{D^{\alpha^i}
\Big(\dfrac{f_j}{f_0}\Big)}{\Big(\dfrac{f_j}{f_0}\Big)}\Bigg|
= o(T_f(r)),\quad 0 \leq i \leq n, 1 \leq j \leq n.
\end{align*}
Hence, we get
\begin{align*}
\Big\Vert \int\limits_{S(r)} \text{log}^+
\Big|\frac{W^\alpha (f_0,\dots,f_n)}{f_0 \cdots f_n}\Big| \sigma =
o(T_f(r)).
\end{align*}
\end{proof}

We finally will need the following estimates of the divisors
of such logarithmic expressions:

\begin{proposition} \label{lemy} {\rm (Special case of \cite{b5},  Proposition 4.10)}
Let $f$ be a linearly nondegenerate meromorphic map of
$\C^m$ into $\C\P^n$ with reduced representation $f=(f_0:\dots : f_n)$. Assume that $\alpha = (\alpha^0, \dots , \alpha^n)$ is an
admissible set for $F=(f_0, \dots , f_n)$, and let again
$p_0 = \min \{ p' : \ell_{F}(p') = n+1\}$.
Then we have
$$\nu_\frac{f_0 \cdot \dots \cdot f_n}{W^{\alpha}(f_0, \dots , f_n)}
\leq \sum_{i=0}^n \min \{\nu_{f_i}, p_0 \} $$ outside an analytic
set of codimension at least two.

\end{proposition}

\section{Regular sequences}
Throughout of this paper, we use the lexicographic order on
$\N_0^p$. Namely, $(i_1,\dots,i_p) > (j_1,\dots,j_p)$
iff for some $s \in \{1,\dots,p\}$ we have $i_\ell = j_\ell$
for $\ell < s$ and $i_s > j_s$.

\begin{lemma}\label{lem2}
Let $A$ be a commutative ring and let
$\{\phi_1,\dots,\phi_p\}$ be a regular sequence in $A$, i.e. for
$i=1,...,p$, $\phi_i$ is not a zero divisor of $A/(\phi_1,...,\phi_{i-1})$.
Denote by $I$ the ideal in $A$ generated by $\phi_1, \dots,\phi_p$.
Suppose that for some $q, q_1, \dots,q_h \in A$ we have an equation
\begin{align*}
\phi_1^{i_1} \cdots \phi_p^{i_p} \cdot q =
\sum_{r=1}^h \phi_1^{j_1(r)} \cdots \phi_p^{j_p(r)}\cdot q_r\:,
\end{align*}
where $(j_1(r),\dots, j_p(r)) > (i_1,\dots,i_p)$ for $r = 1,\dots,h$.
Then $q \in I$.
\end{lemma}
\noindent For the proof, we refer to \cite{Co}, Lemma 2.2. \hfill $\square$ \\

\begin{proposition}\label{reg}

Let $\{Q_j\}_{j=1}^q$ $(q \geq n+1)$ be an admissible set of homogeneous
polynomials of common degree $d \geq 1$ in $\Cal K_f [x_0,\dots,x_n]$.
Then  for any pairwise different
$1 \leq j_0 , \dots , j_n \leq q$ the
sequence $\{ Q_{j_0}, Q_{j_1},..., Q_{j_n} \}$ of elements in
$\Cal K_{\{Q_j\}_{j=1}^q}[x_0,\dots,x_n]$ is a regular sequence,
as well as all its subsequences.
\end{proposition}
\begin{proof}
Since $\Cal K_{\{Q_j\}_{j=1}^q}$ is a field, the ring $\Cal
K_{\{Q_j\}_{j=1}^q}[x_0,\dots,x_n]$ is a local Cohen-Macaulay ring
with maximal ideal $\Cal M = (x_0,...,x_n) \subset \Cal
K_{\{Q_j\}_{j=1}^q}[x_0,\dots,x_n]$ (see for example  \cite{b11},
page 112).  Suppose that $\{ Q_{j_0}, Q_{j_1},..., Q_{j_n} \}$ is a
system of parameters of the ring $\Cal
K_{\{Q_j\}_{j=1}^q}[x_0,\dots,x_n]$, this means (see \cite{b11},
pages 73 and 78) that there exists a natural number $\rho \in \N$
such that
\begin{equation}\label{sp}
{\Cal M}^{\rho} \subset (Q_{j_0}, Q_{j_1},..., Q_{j_n}) \subset
{\Cal M}\;.
\end{equation}
Then by Theorem 31 of \cite{b11}, any subsequence of $\{ Q_{j_0}, Q_{j_1},..., Q_{j_n} \}$ is a regular sequence in $\Cal K_{\{Q_j\}_{j=1}^q}[x_0,\dots,x_n]$.

Since the $\{Q_j\}_{j=1}^q$ $(q \geq n+1)$ are homogeneous
polynomials of common degree $d \geq 1$, the second inclusion
of equation (\ref{sp}) is trivial. In order to prove the first inclusion,
again by Proposition \ref{lang} and its notations
there exists a positive integer $s$
and polynomials $\big\{\widetilde
b_{ik}\big\}_{0 \leq i ,k\leq n}$ in $\Z[T,x]$,
zero or
homogenous in $x$ of degree $s-d$,
such that
\begin{align*}
x_i^s \cdot \widetilde R = \sum_{k=0}^{n} \widetilde b_{ij_k} \widetilde Q_{j_k}\quad
\text{for all}\ i \in \{0,\dots,n\},
\end{align*}
and, since $\big\{Q_{j_k}\big\}_{k=0}^{n}$ is an admissible set,
$R = \widetilde R(\dots, a_{j_kI}, \dots) \not\equiv 0$.
Set
\begin{align*}
b_{ij_k} = \widetilde b_{ij_k}\big((\dots,a_{j_kI},\dots), (x_0,\dots,x_n)\big),\quad
0 \leq i,k \leq n.
\end{align*}
Then it is clear that $R \in \Cal K_{\{Q_j\}_{j=1}^q}$,
$ b_{ij_k} \in  \Cal K_{\{Q_j\}_{j=1}^q}[x_0, \dots , x_n]$.
So we get that
\begin{align*}
x_i^s \cdot  R = \sum_{k=0}^{n}  b_{ij_k}  Q_{j_k}\quad
\text{for all}\ i \in \{0,\dots,n\},
\end{align*}
implying that $x_i^s \in  (Q_{j_0}, Q_{j_1},..., Q_{j_n})$ for all
$i=0,...,n$.  So if take any  $\rho \geq (n+1)(s-1)+1$, then we get
the first inclusion of equation (\ref{sp}), and we are done.
\end{proof}

Let $f$ be a nonconstant meromorphic map of $\C^m$ into $\C P^n$ and
$\big\{Q_j\big\}_{j=1}^q$ $(q \geq n+1)$ be an admissible set of
homogeneous polynomials of degree $d$ in $\Cal K_f[x_0,\dots,x_n]$.
For a nonnegative integer $N$, we denote by $V_N$ the vector space
(over $\Cal K_{\{Q_j\}_{j=1}^q}$)  consisting of all homogeneous
polynomials of degree $N$ in  $\Cal
K_{\{Q_j\}_{j=1}^q}[$ $x_0,\dots,x_n]$ (and of the zero polynomial).
Denote by $(Q_1,\dots,Q_n)$ the ideal in $\Cal K_{\{Q_j\}_{j=1}^q}
[x_0,\dots,x_n]$ generated by $Q_1,\dots,Q_n$.

The following result is similar to Lemma 5 of An-Wang \cite{An}.
However, they proved it for the function field of a smooth projective
variety instead of $\Cal K_{f}$, only for sufficiently big $N$, and with a less elementary method, so we do not try to adopt their proof, but give
an independant one.

\begin{proposition}\label{lem3bis}
Let $\{Q_j\}_{j=1}^q$ $(q \geq n+1)$ be an admissible set of
homogeneous polynomials of common degree $d \geq 1$ in $\Cal K_f
[x_0,\dots,x_n]$. Then for any nonnegative integer $N$ and for any $J:=
\{j_1, \dots , j_n \} \subset \{1, \dots , q\}$, the dimension of
the vector space $\frac{V_N}{(Q_{j_1},\dots,Q_{j_n}) \cap V_N}$ is
 equal to the number of $n$-tuples $(i_1,\dots,i_n)\in
\N_0^n$ such that $i_1+\cdots+i_n\leq N$ and $0\leq
i_1,\dots,i_n\leq d-1.$ In particular, for all $N\geq n(d-1),$ we
have
\begin{align*}
\text{\rm dim}\frac{V_N}{(Q_{j_1},\dots,Q_{j_n}) \cap V_N} = d^n .
\end{align*}

\end{proposition}

\begin{proof} The case $N=0$ holds trivially, so we assume that
$N$ is positive for the rest of the proof.
We first prove that
\begin{align}
\text{\rm dim}\frac{V_N}{(Q_{j_1},\dots,Q_{j_n})\cap V_N}=\text{\rm
dim}\frac{V_N}{(Q_{1},\dots,Q_{n})\cap V_N}\label{e1}
\end{align}
for any choice of $J:=\{j_1,\dots,j_n\}\in\{1,\dots,q\}$ and any
$N.$
 For this it suffices to prove that
$$\text{\rm dim}(Q_1,\dots,Q_n) \cap V_N =
\text{\rm dim}(Q_{j_1},\dots,Q_{j_n}) \cap V_N\; .$$ Since the order
of the $Q_j$ does not matter, it suffices to prove
\begin{equation} \label{e2}
\text{\rm dim}(Q_1,\dots,Q_n) \cap V_N =
\text{\rm dim}(Q_{1},\dots, Q_{n-1}, Q_{j_n}) \cap V_N\: ,
\end{equation}
the rest follows by induction.
But for (\ref{e2}) it suffices to prove:
\begin{equation} \label{e3}
\text{\rm dim}\frac{(Q_1,\dots,Q_n) \cap V_N}
{ (Q_1,\dots,Q_{n-1}) \cap V_N}=
\text{\rm dim}\frac{(Q_{1},\dots, Q_{n-1}, Q_{j_n}) \cap V_N}
{ (Q_1,\dots,Q_{n-1}) \cap V_N}\: .
\end{equation}
We denote for simplicity ${\Cal K}:= {\Cal K}_{\{Q_j\}_{j=1}^q}$ and
let $\phi$ be the following ${\Cal K}$-linear map:
$$\phi :
\frac{(Q_1,\dots,Q_n) \cap V_N}
{ (Q_1,\dots,Q_{n-1}) \cap V_N} \rightarrow
\frac{(Q_{1},\dots, Q_{n-1}, Q_{j_n}) \cap V_N}
{ (Q_1,\dots,Q_{n-1}) \cap V_N}\; ; $$ $$
[\sum_{j=1}^{n-1}b_jQ_j +b_nQ_n ] \mapsto
[\sum_{j=1}^{n-1}b_jQ_j +b_nQ_{j_n} ]
$$
with $b_j \in {\Cal K}[x_0, \dots , x_n]$.
This map is clearly surjective, so if we still prove that it is
well defined and injective, we get (\ref{e3}).
In order to prove that $\phi$ is well defined, let
$[\sum_{j=1}^{n-1}b_jQ_j +b_nQ_n ]=[\sum_{j=1}^{n-1}b_j'Q_j +b_n'Q_n ]$. This means that
$(b_n-b_n')Q_n \in (Q_1,\dots,Q_{n-1}) \cap V_N$.
But since by Proposition \ref{reg}, $Q_1, \dots , Q_n$ is
a regular sequence, $Q_n$ is not a zero divisor in
$\frac{{\Cal K}[x_0, \dots , x_n]}{(Q_1, \dots , Q_{n-1})}$,
so that $(b_n - b_n') \in (Q_1, \dots , Q_{n-1})$. Hence,
$$[\sum_{j=1}^{n-1}b_jQ_j +b_nQ_{j_n} ]-
[\sum_{j=1}^{n-1}b_j'Q_j +b_n'Q_{j_n }] = [(b_n - b_n')Q_{j_n}] =
0$$ in $\frac{(Q_{1},\dots, Q_{n-1}, Q_{j_n}) \cap V_N} {
(Q_1,\dots,Q_{n-1}) \cap V_N}$, so $\phi$ is well defined. The
injectivity of $\phi$ follows by the same argument, just changing
the roles of $Q_n$ and $Q_{j_n}$, since by Proposition~\ref{reg},
$Q_1, \dots , Q_{n-1}, Q_{j_n}$ is also a regular sequence. Hence,
we get (\ref{e2}) and, thus, (\ref{e1}). We finally remark that for the proof of (\ref{e2}) we only used that $\{Q_1,...Q_n, Q_{j_n}\}$ is an 
admissible set of homogenous polynomials of common degree $d$.

Take a point $z_0\in\C^m$ such that the hypersufaces in $\C P^n$
defined by $Q_1(z_0),\dots,$ $Q_{n+1}(z_0)$ have no common point.
Since $Q_1(z_0),\dots,Q_n(z_0)$ define a subvariety of dimension 0,
there exists a hyperplane $H_{n+1}$ in $\C P^n$ such that
$\cap_{i=1}^n Q_j(z_0)\cap H_{n+1}=\varnothing.$ Furthermore, by
induction, there exist hyperplanes $H_1,\dots,H_{n+1}$ such that
$\cap_{j=1}^{i-1}Q_j(z_0)\cap_{k=i}^{n+1}H_k=\varnothing,$ for all
$i\in\{1,\dots,n+1\}.$ This means that $\{Q_1,\dots,Q_{i-1},
H_{i}^d,\dots, H_{n+1}^d\}$ is an admissible set, for any
$i\in\{1,\dots, n+1\}.$ Then, by (\ref{e2}), taking into accont the remark at the end of its proof,  and by induction, we get
that
\begin{align}
\text{\rm dim}\frac{V_N}{(Q_{1},\dots,Q_{n}) \cap V_N} = \text{\rm
dim}\frac{V_N}{(H_1^d,\dots,H_n^d) \cap V_N}.\label{new1}
\end{align}
As  $H_1,\dots, H_{n+1}$ are linearly independent, it follows from a
well-known fact of linear algebra that there exists a permutation
$\{k_1,\dots,k_{n+1}\}$ of $\{0,\dots,n\}$ such that
$H_1,\dots,H_{i-1}, x_{k_i},\dots,x_{k_{n+1}}$ are linearly
independent, for any $i\in\{1,\dots,n+2\}.$ This means that
$\{H_1^d,\dots,H_{i-1}^d, x_{k_i}^d,\dots,x_{k_{n+1}}^d\}$ is an
admissible set. Then, by (\ref{e2}) and by induction, we get that
\begin{align}
 \text{\rm
dim}\frac{V_N}{(H_1^d,\dots,H_n^d) \cap V_N}=\text{\rm
dim}\frac{V_N}{(x_1^d,\dots,x_n^d) \cap V_N}.\label{new0}
\end{align}
By (\ref{e1}), (\ref{new1}) and (\ref{new0}), for all positive
integer $ N$ we have
\begin{align*}
\text{\rm dim}\frac{V_N}{(Q_{j_1},\dots,Q_{j_n}) \cap V_N}
=\text{\rm dim}\frac{V_N}{(x_1^d,\dots,x_n^d) \cap V_N}.
\end{align*}
On the other hand, it is easy to see that for any positive integer
$N,$ the vector space $\frac{V_N}{(x_1^d,\dots,x_n^d) \cap V_N}$ has
a basis $\{[x_0^{N-(i_1+\cdots i_n)}x_1^{i_1}\cdots x_n^{i_n}],
i_1+\cdots+i_n\leq N, 0\leq i_1,\dots, i_n\leq d-1\}.$ This
completes the proof of Proposition \ref{lem3bis}.
\end{proof}

\section{Proof of Main Theorem}

We first  prove the theorem for the case where all the $Q_j$
$(j=1,\dots,q)$ have the same degree $d$.

We may assume, without loss of generality, that $f$ is algebraically
nondegenerate over ${\Cal K}_{\{Q_j\}_{j=1}^q}$: We replace the
polynomials $\big\{Q_j\big\}_{j=1}^q$ by the polynomials $\tilde
Q_{j} :=\frac{1}{a_{jI_{2}}}Q_j$, where $a_{jI_{2}}\not=0$ is any
nonzero coefficient of $Q_{j}$, $j=1,...,q$. Then
 $\big\{\tilde Q_j\big\}_{j=1}^q$ is also an admissible set of homogeneous
polynomials in $\Cal K_f [x_0,\dots,x_n]$ with $\text{deg}\,Q_j = d \geq 1$.
Since $f$ is algebraically nondegenerate over
$\tilde {\Cal K}_{\{Q_j\}_{j=1}^q}$
and
$\tilde {\Cal K}_{\{Q_j\}_{j=1}^q} \supset
{\Cal K}_{\{\tilde Q_j\}_{j=1}^q}$,
we have that $f$ is algebraically nondegenerate over
${\Cal K}_{\{\tilde Q_j\}_{j=1}^q}$.
So we get that
 for any $\varepsilon > 0$, there exists a positive integer $L$, depending on $n$, $\epsilon$ and $d$ in an explicit way, such that
\begin{align*}
\Vert (q-n-1-\varepsilon) T_f(r) \leq \sum_{j=1}^q \frac{1}{d}
N^{(L)}_f(r,\tilde Q_j).
\end{align*}
But since
\begin{align*}
N^{(L)}_f(r,\tilde Q_j) &= N^{(L)}_f(r,\frac{1}{a_{jI_{2}}}Q_j)
=
N^{(L)}_{\frac{1}{a_{jI_{2}}}Q_j \circ f}(r)\\
&\leq  N_{\frac{1}{a_{jI_{2}}}}(r) + N^{(L)}_{Q_j \circ f}(r) =
N^{(L)}(r, Q_j \circ f) + o(T_{f}(r))\, ,
\end{align*}
we finally get that
for any $\varepsilon > 0$ we have
\begin{align*}
\Vert (q-n-1-\varepsilon) T_f(r) \leq \sum_{j=1}^q \frac{1}{d}
N^{(L)}_f(r, Q_j).
\end{align*}\\

For each nonnegative integer $k$, we denote again by $V_k$ the space
(over $\Cal K_{\{Q_j\}_{j=1}^q}$) of homogeneous polynomials
of degree $k$  (and of the zero polynomial) in $\Cal K_{\{Q_j\}_{j=1}^q}[x_0,$ $\dots,x_n]$.
Set $V_k = \{0\}$ for $k < 0$.

Let $J:= \{j_1, \dots , j_n \} \subset \{1, \dots , q\}$.
For each positive integer $N$ divisible by $d$ and for each
$I := (i_1,\dots,i_n) \in \N_0^n$ with
$\Vert I \Vert := \sum\limits_{s=1}^n i_s \leq \dfrac{N}{d}$,
we set
 \begin{align*}
V_N^I = \sum_{E:= (e_1,\dots,e_n) \geq I}
Q_{j_1}^{e_1} \cdots Q_{j_n}^{e_n} \cdot V_{N-d\Vert E\Vert}.
\end{align*}
Note that $V_N^I \supset V_N^J$ if $I < J$ (lexicographic order),
and $V_N^{(0,\dots,0)} = V_N$.

Denote by $\{I_1,\dots,I_K\}$ the set of all $I \in \N_0^n$ with
$\Vert I \Vert \leq \dfrac{N}{d}$. We write $I_k =
(i_{1k},\dots,i_{nk})$, \ $k = 1,\dots,K$. Assume that $I_1 =
(0,\dots,0) < I_2 < \dots < I_K=(\frac{N}{d},0,...,0)$. We have
\begin{align*}
V_N = V_N^{I_1} \supset V_N^{I_2} \supset \dots \supset V_N^{I_K}
= \{u \cdot Q_{j_1}^{\frac{N}{d}}: u \in \Cal K_{\{Q_j\}_{j=1}^q}\} \\ \text{and}\ K = K(N,d,n) = \begin{pmatrix} \dfrac{N}{d}+n\cr
n \end{pmatrix}.
\end{align*}
Set
\begin{align*}
m_k := \text{dim} \frac{V_N^{I_k}}
{V_N^{I_{k+1}}},\quad k = 1, \dots, K-1,\ \text{and}\ m_K:= \text{dim}
V_N^{I_K}=1.
\end{align*}

We now prove that:\\
Although the $V_N^{I_k}$ may depend on $J$, the $m_k$, $k=1, \dots ,
K$ are independent of $J$. Moreover,
\begin{align}
m_k = d^n \label{eq31}
\end{align}
for all $N$ divisible by\ $d$ and for all $k \in \{1, \dots, K\}$
with $N - d\Vert I_k\Vert \geq nd.$

We define vector space homomorphisms
\begin{align*}
\varphi_k : V_{N-d\Vert I_k\Vert} \longrightarrow
\frac{V_N^{I_k}}{V_N^{I_{k+1}}}\quad (k=1,\dots, K-1)
\end{align*}
as $\varphi_k(\gamma) = \big[Q_{j_1}^{i_{1}k} \cdots Q_{j_n}^{i_{n}k} \gamma \big]$,
where $\gamma \in V_{N-d\Vert I_k\Vert}$ and
$\big[ Q_{j_1}^{i_{1}k} \cdots Q_{j_n}^{i_{n}k}\gamma\big]$
is the class in $\dfrac{V_N^{I_k}}{V_N^{I_{k+1}}}$ containing
$Q_{j_1}^{i_1k} \cdots Q_{j_n}^{i_nk}\gamma$.

It is clear that the $\varphi_k$ $(1 \leq k \leq K-1)$ are surjective (note that
for any $E \in \N_0^n$ with $\Vert E \Vert \leq \dfrac{N}{d}$ and
$E > I_k$ then $E \geq I_{k+1}$).

For any $\gamma \in \text{ker}\,\varphi_k$
\begin{align*}
Q_{j_1}^{i_1k} \cdots Q_{j_n}^{i_nk}\gamma
&\in \sum_{E=(e_1,\dots,e_n) \geq I_{k+1}}
Q_{j_1}^{e_1} \cdots Q_{j_n}^{e_n} V_{N-d\Vert E\Vert}\\
&= \sum_{E=(e_1,\dots,e_n) > I_k}
Q_{j_1}^{e_1} \cdots Q_{j_n}^{e_n} V_{N-d\Vert E\Vert}.
\end{align*}
So we have
\begin{align*}
Q_{j_1}^{i_1k} \cdots Q_{j_n}^{i_nk}\gamma =
\sum_{E=(e_1,\dots,e_n) > I_k}
Q_{j_1}^{e_1} \cdots Q_{j_n}^{e_n} \gamma_E\:,
\end{align*}
where $\gamma_E \in V_{N-d\Vert E\Vert}$. Furthermore,
 by Lemma \ref{lem2} and Proposition \ref{reg} we have $\gamma \in (Q_{j_1},\dots,Q_{j_n})$. Thus

\begin{align}
\text{ker}\,\varphi_k \subset (Q_{j_1},\dots,Q_{j_n}) \cap V_{N-d\Vert I_k\Vert}\:.
\label{eq32}
\end{align}
Conversely, for any $\gamma \in (Q_{j_1},\dots,Q_{j_n}) \cap V_{N-d\Vert I_k\Vert}$
$(\gamma \neq 0)$,
\begin{align*}
\gamma = \sum_{s=1}^n \gamma_s Q_{j_s},\quad
\gamma_s \in V_{N-d(\Vert I_k\Vert +1)},
\end{align*}
we have
\begin{align*}
 I_s' := (i_{1k},\dots,i_{sk}+1, \dots,i_{nk}) > I_k,\quad
(s = 1,\dots,n)
\end{align*}
and $\Vert  I_s'\Vert = \Vert I_k\Vert +1 \leq \dfrac{N}{d}$
(since $ \gamma \neq 0$).
So we get $ I_s' \geq I_{k+1}$, $s = 1, \dots,n$. Thus
\begin{align*}
Q_{j_1}^{i_{1k}}\cdots Q_{j_n}^{i_{nk}} \gamma =
\sum_{s=1}^n Q_{j_1}^{i_{1k}} \cdots Q_{j_s}^{i_{sk}+1} \cdots Q_{j_n}^{i_{nk}}
\cdot \gamma_s \in V_N^{I_{k+1}}\:.
\end{align*}
This means that $\gamma \in \text{ker}\,\varphi_k$. So we have
\begin{align}
\text{ker}\,\varphi_k \supset (Q_{j_1},\dots,Q_{j_n}) \cap V_{N-d \Vert I_k\Vert}\:.
\label{eq33}
\end{align}

By (\ref{eq32}), (\ref{eq33}) and since $\varphi_k$ is surjective, we have:

\begin{equation} \label{ll}
m_k = \text{dim} \frac{V_N^{I_k}}{V_N^{I_{k+1}}} \simeq
\text{dim} \frac{V_{N-d\Vert I_k\Vert}}
{(Q_{j_1},\dots,Q_{j_n}) \cap V_{N-d\Vert I_k \Vert}}\: ,\:k \in \{1,...,K-1 \}.
\end{equation}
Hence, by Proposition \ref{lem3bis} we get (\ref{eq31}) and the
independence of $J$ of $m_k.$ \hfill $\square$

Since $V_N = V_N^{I_1} \supset V_N^{I_2} \supset \dots \supset
V_N^{I_K} $ and $m_k = \text{dim} \dfrac{V_N^{I_k}}
{V_N^{I_{k+1}}}, (k \in \{1,...,K-1\})$, $ m_K=\text{dim} V_N^{I_K}=1$, we may choose a basis
$\{\psi_1^J,\dots,\psi_M^J\}$ ($M = \begin{pmatrix} N+n\cr
n\end{pmatrix}$) of $V_N$ such that
\begin{align*}
\Big\{ \psi_{M-(m_k+\dots+m_{K})+1}^J,\dots,\psi_M^J\Big\}
\end{align*}
is a basis of $V_N^{I_k}$ for any $k \in \{1,\dots,K\}$. For each
$k \in \{1,\dots,K\}$ and $ \ell \in \big\{ M - (m_{k+1} + \dots +
m_{K}),\dots, M - (m_k+\dots+m_{K})+1\big\},$ we have
\begin{align}
\psi_\ell^J = Q_{j_1}^{i_{1k}} \cdots Q_{j_n}^{i_{nk}}
\gamma_\ell^J,\; \text{where}\; \gamma_\ell^J \in V_{N-d\Vert I_k
\Vert}.\label{new10}
\end{align}
\noindent Then, we have
\begin{align}
\prod_{j=1}^M \psi_j^J(f) = \prod_{k=1}^{K}
\big((Q_{j_1}(f))^{i_{1k}} \cdots (Q_{j_n}(f))^{i_{nk}}\big)^{m_k} \cdot\prod_{\ell=1}^M\gamma_\ell^J(f)\label{44}
\end{align}
By Lemma \ref{lem1} there exists $h_\ell^J \in \Cal C_f$ such that,
outside an analytic subset in $\C^m$ of codimension at least two,
\begin{align*}
|\gamma_\ell^J (f)| \leq h_\ell^J\cdot \Vert f\Vert^{N-d\Vert I_k\Vert}.
\end{align*}
So we get
\begin{align*}
\prod_{\ell=1}^M |\gamma_\ell^J(f)| \leq \prod_{k=1}^{K}
\big(
\Vert f\Vert^{N-d\Vert I_k\Vert}\big)^{m_k} \cdot h^J\:,
\end{align*}
where $h^J:=\prod_{\ell=1}^Mh_\ell ^J\in \Cal C_f$. This implies
that (outside a proper analytic subset of $\C^m$)
\begin{align}
\text{log} \prod_{\ell=1}^M |\gamma_\ell^J(f)|
\leq \sum_{k=1}^K m_k(N-d\Vert I_k\Vert)
\text{log}\Vert f\Vert + \text{log}h^J. \label{eq34}
\end{align}

By (\ref{ll}) and since $m_K=1$, we have that $m_k$ only depends on $\Vert I_k \Vert,$  i.e. $m_k=m(\Vert I_k \Vert)$,
$k=1, \dots , K$. So we have, for $s=1,...,n$,
\begin{align*}
\sum_{k=1}^{K}m_k\cdot i_{sk} =
\sum_{\ell=0}^{\frac{N}{d}}\sum_{k:\Vert I_k \Vert = \ell}m_k\cdot i_{sk} =\sum_{\ell=0}^{\frac{N}{d}}m(\ell)\sum_{k:\Vert I_{k}\Vert=\ell}i_{sk}\,.
\end{align*}
Now for every $\ell$ the  the symmetry $(i_1,\cdots, i_n)\to
((i_{\sigma(1)}, \dots,i_{\sigma(n)})$ shows that $\sum_{k:\Vert
I_{k}\Vert=\ell}i_{sk}$ is independent of $s$. So, we get
\begin{align}
A:=\sum_{k=1}^{K}m_k\cdot
i_{sk}\quad \text{is independent of}\quad s \quad \text{and} \quad
J\, ,\label{51}
\end{align}
the latter by (\ref{eq31}).

Denote by $\mathcal B$ the set of all $k\in\{1,\dots,K\}$ such that
$N-d\Vert I_k\Vert \geq nd$.

\noindent Put
\begin{align*}
\widetilde I_k := (i_{1k}, \dots,i_{nk}, i_{(n+1)k}),\quad k
\in\mathcal B,
\end{align*}
where $i_{(n+1)k} := (\dfrac{N}{d}-n) - (i_{1k}+\dots+i_{nk})$. Then
$\{\widetilde I_k:\; k\in\mathcal B\}$ is the set of all $\widetilde
I \in \N_0^{n+1}$ with $\Vert \widetilde I\Vert = \dfrac{N}{d}-n$.
For any $\widetilde I := (i_1,\dots,i_{n+1}) \in \{\widetilde I_k:\;
k\in\mathcal B\}$  and for any bijection $\sigma : \{1,\dots,n+1\}
\to \{1,\dots,n+1\}$, we have $(i_{\sigma(1)},
\dots,i_{\sigma(n+1)}) \in \{\widetilde I_k:\; k\in\mathcal B\}.$
Therefore, by (\ref{eq31}) we have
\begin{align}
A\geq\sum_{k \in \mathcal B} m_k\cdot i_{sk}=d^n\sum_{k \in \mathcal
B}  i_{sk}=\frac{d^n}{n+1}\sum_{k \in \mathcal B}\Vert
\widetilde I_k\Vert=\frac{d^n}{n+1}\begin{pmatrix} \dfrac{N}{d}\cr n
\end{pmatrix}(\frac{N}{d}-n)
 .\label{51a}
\end{align}
We have
\begin{align}
\sum_{k=1}^{K}m_k\cdot (\frac{N}{d}- \Vert I_k\Vert )=
\sum_{k=1}^{K}m_k\cdot \frac{N}{d}-\sum_{k=1}^{K}m_k\cdot \Vert
I_k\Vert = \frac{MN}{d}-nA. \label{54}
\end{align}
By (\ref{44}), we have  for $N$ divisible by $d$:
\begin{align}
\prod_{j=1}^M \psi_j^J(f) =\big(Q_{j_1}(f) \cdots Q_{j_n}(f)\big)^A
\cdot\prod_{l=1}^M\gamma_\ell^J(f)\label{45}
\end{align}
By  (\ref{eq34}), (\ref{54}) and (\ref{45}), we have (outside a
proper analytic subset of $\C^m$) for $N$ divisible by $d$ :
\begin{align*}
\text{log}\prod_{j=1}^M |\psi_j^J(f)| \leq A\cdot
\text{log}\prod_{i=1}^n |Q_{j_i}(f)|
+(\frac{MN}{d}-nA)d\cdot\text{log}\Vert f \Vert + \text{log}h^J\,.
\end{align*}
If we still choose the function $h \in \Cal C_f$, with $h \geq 1$, common for all $J$,
for example by putting $h:=\prod_J (1+h_J)$),
 we get
\begin{align}
\text{log}\prod_{i=1}^n |Q_{j_i}(f)| \geq \frac{1}{A}
(\text{log}\prod_{j=1}^M |\psi_j^J(f)| - \text{log}h)
-(\frac{MN}{dA}-n) d \cdot \text{log}\Vert f \Vert.\label{eq36}
\end{align}
We choose $N:=d\cdot[2(n+1)(2^n-1)(nd+1)\epsilon^{-1}+n+1].$ Then
 by (\ref{51a}), we have (assuming without loss of generality that 
 $\epsilon < 1$)
 \begin{align}
d \cdot (\frac{MN}{dA}-n-1)\leq d\cdot\big(\frac{N\begin{pmatrix}
N+n\cr n
\end{pmatrix}}{\frac{d^{n+1}}{n+1}\begin{pmatrix} \dfrac{N}{d}\cr n
\end{pmatrix}(\frac{N}{d}-n)}-n-1\big)\notag\\
=d(n+1)\cdot\big(\prod_{i=1}^n\frac{N+i}{N-(n+1-i)d}-1\big)<d(n+1)\big((\frac{N+1}{N-nd})^n-1\big)\notag\\
=d(n+1)\big((1+\frac{nd+1}{N-nd})^n-1\big)<d(n+1)(2^n-1)\frac{nd+1}{N-nd}\notag\\
\leq
d(n+1)(2^n-1)\frac{nd+1}{d\cdot\big(2(n+1)(2^n-1)(nd+1)\epsilon^{-1}+n\big)-nd}=\frac{\epsilon}{2}.\label{new2}
 \end{align}
 By (\ref{eq36}) and Lemma \ref{lem1} (applied to every factor
$Q_{\beta_j}$, $j=1, \dots , q-n$, using that we  can complete every $Q_{\beta_j}$ with $n$ other $Q_j$ not having bigger norm, so that the maximum of the norms is obtained by $Q_{\beta_j}$), we have
\begin{align*}
&\text{log}\prod_{j=1}^q |Q_j(f)| =
\max_{\{\beta_1,\dots,\beta_{q-n}\} \subset \{1,\dots,q\}}
\text{log}|Q_{\beta_1}(f) \cdots Q_{\beta_{q-n}}(f)| \\
&\qquad + \min_{J = \{j_1,\dots,j_n\} \subset \{1,\dots,q\}}
\text{log}|Q_{j_1}(f) \cdots Q_{j_n}(f)|\\
&\geq (q-n)d \cdot \text{log}\Vert f\Vert + \min_{J \subset
\{1,\dots,q\}} \frac{1} {A}
\text{log}\prod_{j=1}^M |\psi_j^J(f)|\\
&\qquad - d \cdot(\frac{MN}{Ad}-n) \text{log}\Vert f\Vert
 - \text{log}\widetilde h\\
&= (q-n-1)d \cdot \text{log}\Vert f \Vert + \frac{1} {A}
\min_{J \subset \{1,\dots,q\}}\text{log} \prod_{j=1}^M |\psi_j^J(f)|\\
&\qquad - d \cdot (\frac{MN}{Ad}-n-1) \cdot \text{log}\Vert f \Vert
- \text{log}\widetilde h,
\end{align*}
where  the choices of the indices for the maximum respectively the
minimum may depend on $z$, however, by (observing $A \geq 1$ and by) choosing $\widetilde h $ as a
product of the form $\prod (1+h_{\nu})$, where the $h_{\nu}$ run
over all the possible choices, we obtain $\widetilde h \in \Cal
C_f$. Furthermore we observe that the first and the last term are
well defined outside an analytic subset of $\C^m$ of codimension at
least two and the choices of maxima and minima are locally finite
there, in particular the resulting functions are continuous there as
functions with values in $[0,+\infty]$. Hence, the inequality still
holds outside an analytic subset of $\C^m$ of codimension at least
two  by continuity. So by integrating and by using (\ref{new2}),
outside an analytic subset of codimension at least two in $\C^m
\supset S(r)$ we get
\begin{align}
\int\limits_{S(r)} \text{log} \prod_{j=1}^q |Q_j(f)| \sigma \geq
(q-n-1) d \cdot T_f(r) &+ \frac{1}{A} \int\limits_{S(r)} \min_J
\text{log} \prod_{j=1}^M
|\psi_j^J(f)| \sigma \notag\\
& - \frac{\varepsilon}{2}T_f(r) - o(T_f(r)).\label{eq37}
\end{align}

We write
\begin{align*}
\psi_j^J = \sum_{I \in \Cal T_N} c_{jI}^{J}x^I \in V_N,\quad
c_{jI}^J \in \Cal K_{\{Q_j\}_{j=1}^q}
\end{align*}
with $j=1,\dots,M$, \ $J \subset \{1,\dots,q\}$, \ $\#J = n$. For
each $j \in \{1,\dots,M\}$ and $J \subset \{1,\dots,q\}$, $\# J = n$
we fix an index $I_{i}^J \in \Cal T_N$ such that $c_{jI_{j}^J}^J
\not\equiv 0$. Define
\begin{align*}
\xi_{jI}^J = \frac{c_{jI}^J}{c_{jI_{j}^J}^J},\quad j \in
\{1,\dots,M\}, \ J \subset \{1,\dots,q\}, \ \#J = n.
\end{align*}
Set
\begin{align*}
\widetilde \psi_j^J = \sum_{I \in \Cal T_N} \xi_{jI}^Jx^I
\in \Cal K_{\{Q_j\}_{j=1}^q} [x_0,\dots,x_n].
\end{align*}
For each positive integer $p$, we denote by $\Cal L(p)$ the vector
space generated over $\C$ by
\begin{align*}
\Big\{ \prod_{\begin{matrix}
\scriptstyle{1 \leq j \leq M, I \in \Cal T_N}\cr
\noalign{\vskip-0.15cm}
\scriptstyle{J \subset \{1,\dots,q\}, \#J =n}\end{matrix}}
\big(\xi_{jI}^J\big)^{n_{jI}^J} : n_{jI}^J \in \N_0,
\sum_{\begin{matrix}
\scriptstyle{1 \leq j \leq M, I \in \Cal T_N}\cr
\noalign{\vskip-0.15cm}
\scriptstyle{J \subset \{1,\dots,q\}, \# J = n}\end{matrix}}
n_{jI}^J = p \Big\}\:.
\end{align*}
We have \ $\Cal L(p) \subset \Cal L(p+1) \subset \Cal
K_{\{Q_j\}_{j=1}^q}$ (note that $\xi_{jI_{j}^J}^J \equiv 1$, $j \in
\{1,\dots,M\}$, $J \subset \{1,\dots,q\}$, $\#J =n$). Let
$\{b_1,\dots, b_{t_{p+1}}\}$ be a basis of $\Cal L(p+1)$ such that
$\{b_1, \dots,b_{t_p}\}$ is a basis  of $\Cal L(p)$.
 It is easy to see that
 \begin{align}
 t_{p+1}\leq \left(\begin{array}{c}
\binom{n+N}{n}^2.\binom{q}{n} + p\\
\binom{n+N}{n}^2.\binom{q}{n}-1
\end{array}\right)
<\Bigg(\binom{n+N}{n}^2.\binom{q}{n}+p\Bigg)^{
\binom{n+N}{n}^2.\binom{q}{n}-1}\label{new4}
 \end{align} for all positive
integer $p$ (note that $\#\Cal T_N =\binom{n+N}{n}=M$).

 Since
$\big\{ \widetilde \psi_1^J, \dots, \widetilde \psi_M^J\big\}$ is a
basis of $V_N$, $\{b_1,\dots,b_{t_{p+1}}\}$ is a basis of $\Cal
L(p+1)$ and $f$ is algebraically nondegenerate over $\Cal
K_{\{Q_j\}_{j=1}^q}$ we have that $b_k \widetilde \psi_j^J(f)$ ($1
\leq k \leq t_{p+1}$, $1 \leq j \leq M$) are linearly independent
over $\C$. It is easy to see that $b_\ell \widetilde \psi_j^J(f)$
($1 \leq \ell \leq t_p$, $1 \leq j \leq M$, $J \subset
\{1,\dots,q\}$, $\# J = n$) are linear combinations of $b_kf^I$ ($1
\leq k \leq t_{p+1}$, $I \in \Cal T_N$) over $\C$. So for each $J
\subset  \{1,\dots,q\},\; \# J = n$ there exists $A_J \in
\text{mat}(t_{p+1}M \times t_pM, \C)$ such that
\begin{align*}
\big(b_\ell  \widetilde \psi_j^J(f), 1 \leq \ell \leq t_p, 1 \leq j
\leq M\big) = \big(b_kf^I, 1 \leq k \leq t_{p+1}, I \in \Cal T_N)
\cdot A_J
\end{align*}
(note that $\# \Cal T_N = M$). Since $b_\ell \widetilde \psi_j^J(f)$
$(1 \leq \ell \leq t_p, \ 1 \leq j \leq M)$ are linearly independent
over $\C$, we obtain rank\,$A_J = t_pM$.

\noindent Take matrices $B_J \in \text{mat}(t_{p+1}M \times
(t_{p+1}-t_p)M,\C)$ $(J \subset \{1,\dots,q\}, \# J=n)$ such that
\begin{align*}
C_J := A_J \cup B_J  \in GL (t_{p+1}M,\C).
\end{align*}
We write
\begin{align}
\big(b_kf^I , 1 \leq k \leq t_{p+1}, I \in \Cal T_N\big) \cdot C_J
&= \big(b_\ell \widetilde \psi_j^J(f),
1 \leq \ell \leq t_p, \ 1 \leq j \leq M\notag\\
& \qquad h^J_{uv}, t_p+1  \leq u \leq t_{p+1},\ 1 \leq v \leq M
\big)\,. \label{eq38}
\end{align}
Since $\{b_1,\dots,b_{t_{p+1}}\}$ is a basis of $\Cal L(p+1)$ and
$f$ is algebraically nondegenerate over $\Cal K_{\{Q_j\}_{j=1}^q}$
we have that $b_kf^I$ $(1 \leq k \leq t_{p+1}, I \in \Cal T_N)$ are
linearly independent over $\C$. By Lemma \ref{lem5} there exists an
admissible set $\alpha := (\alpha^0,\dots,\alpha^{t_{p+1}M})$  for
$(b_kf^I, 1 \leq k \leq t_{p+1}, I \in \Cal T_N)$. By (\ref{eq38})
we have that $\alpha$ is also an admissible set for
\begin{align*}
\big(b_\ell \widetilde \psi_j^J(f), 1 \leq \ell \leq t_p, 1 \leq j
\leq M, h_{uv}^J , t_p+1 \leq u \leq t_{p+1}, 1 \leq v \leq
M\big)\,.
\end{align*}
Set $W^\alpha := W^\alpha (b_kf^I, 1 \leq k \leq t_{p+1}, I \in \Cal
T_N)$ and $$ W^\alpha_J := W^\alpha(b_\ell \widetilde \psi_j^J(f), 1
\leq \ell \leq t_p, 1 \leq j \leq M, h_{uv}^J, t_p+1 \leq u \leq
t_{p+1}, 1 \leq v \leq M).$$ We have \ $W_J^\alpha = \text{det} C_J
\cdot W^\alpha$.

By Lemma \ref{lem1} we have (with the same arguments on domains
of definition and continuity as above)
\begin{align*}
\int\limits_{S(r)} \min_J \text{log} \prod_{j=1}^M
|\psi_j^J(f)|^{t_p} \sigma &\geq \int\limits_{S(r)} \min_J
\text{log} \prod_{\begin{matrix} \scriptstyle{1 \leq j \leq M}\cr
\noalign{\vskip-0.1cm} \scriptstyle{1 \leq \ell \leq t_p}
\end{matrix}}|b_\ell \widetilde \psi_j^J(f)|\sigma\\
&\quad + \int\limits_{S(r)} \min_J \text{log}
\frac{\prod\limits_{j=1}^M |c_{jI_{j}^J}^J|^{t_p}}
{\prod\limits_{\ell=1}^{t_p} |b_\ell|^M} \sigma
\end{align*}
 \begin{align*}
&\geq  \int\limits_{S(r)} \min_J \text{log} \prod_{\begin{matrix}
\scriptstyle{1 \leq j \leq M}\cr \noalign{\vskip-0.15cm}
\scriptstyle{1 \leq \ell \leq t_p}\end{matrix}}
|b_\ell \widetilde \psi_j^J (f)| \sigma - o(T_f(r))\\
&\geq \int\limits_{S(r)} \min_J \text{log}
\Big(\prod_{\begin{matrix} \scriptstyle{1 \leq j \leq M}\cr
\noalign{\vskip-0.15cm} \scriptstyle{1 \leq \ell \leq
t_p}\end{matrix}}|b_\ell \widetilde \psi_j^J(f)| \cdot
\prod_{\begin{matrix} \scriptstyle{t_p+1\leq u \leq t_{p+1}}\cr
\noalign{\vskip-0.15cm}
\scriptstyle{1 \leq v \leq M}\end{matrix}}|h_{uv}^J|\Big)\sigma\\
&\quad - \int\limits_{S(r)} \max_J \text{log} \prod_{\begin{matrix}
\scriptstyle{t_p+1\leq u \leq t_{p+1}}\cr \noalign{\vskip-0.15cm}
\scriptstyle{1 \leq v \leq M} \end{matrix}}|h_{uv}^J | \sigma -
o(T_f(r))
\end{align*}
\begin{align}
&\geq \int\limits_{S(r)} \min_J \text{log} \Big(
\prod_{\begin{matrix} \scriptstyle{1 \leq j \leq M}\cr
\noalign{\vskip-0.15cm} \scriptstyle{1 \leq \ell \leq
t_p}\end{matrix}}
 |b_\ell \widetilde \psi_j^J (f)|
\cdot \prod_{\begin{matrix} \scriptstyle{t_p+1 \leq u \leq
t_{p+1}}\cr \noalign{\vskip-0.15cm}
\scriptstyle{1 \leq v \leq M}\end{matrix}} |h_{uv}^J|\Big)\sigma \notag\\
&\quad - \int\limits_{S(r)} \text{log}\Vert f
\Vert^{NM(t_{p+1}-t_p)}\sigma -
o(T_f(r))\qquad\qquad\notag\\
&\geq \int\limits_{S(r)} \min_J \text{log}\Big(
\prod_{\begin{matrix} \scriptstyle{1\leq j \leq M}\cr
\noalign{\vskip-0.15cm} \scriptstyle{1\leq \ell \leq
t_p}\end{matrix}}
 |b_\ell \widetilde \psi_j^J(f)| \cdot
\prod_{\begin{matrix} \scriptstyle{t_p+1\leq u \leq t_{p+1}}\cr
\noalign{\vskip-0.15cm}
\scriptstyle{1 \leq v \leq M}\end{matrix}} |h_{uv}^J|\Big)\sigma \notag\\
&\quad - NM(t_{p+1}-t_p)T_f(r) - o(T_f(r))\,,\label{eq39}
\end{align}
where $\min\limits_J$ is taken over  all subset $J \subset \{1,\dots,q\}$,
$\# J = n$.

We may choose a positive integer $p$ such that
\begin{align} p\leq p_0:=\big[\frac{\big( \binom{n+N}{n}^2.\binom{q}{n}-1
\big).\log\big(
\binom{n+N}{n}^2.\binom{q}{n}\big)}{\log(1+\frac{\epsilon}{2MN})}+1\big]^2
\;\text{and}\; \frac{t_{p+1}}{t_{p}}
<1+\frac{\epsilon}{2MN}.\label{new3}
\end{align}
 Indeed, otherwise $\frac{t_{p+1}}{t_{p}}
\geq 1+\frac{\epsilon}{2MN}$ for all $p \leq p_0$. This implies that
$t_{p_0+1} \geq (1+\frac{\epsilon}{2MN})^{p_0}$.
 Therefore, by (\ref{new4}) we have
\begin{align*}
\log(1+\frac{\epsilon}{2MN}) \leq \dfrac{\log t_{p_0+1}}{p_0} &<
\dfrac{\big(\binom{n+N}{n}^2.\binom{q}{n}-1 \big).\log\big(
\binom{n+N}{n}^2.\binom{q}{n}+ p_0\big)}{p_0}\notag\\
&\leq \dfrac{\big( \binom{n+N}{n}^2.\binom{q}{n}-1 \big).\log\big(
\binom{n+N}{n}^2.\binom{q}{n}\big)\cdot\log p_0}{p_0}\notag\\
&\leq \dfrac{\big( \binom{n+N}{n}^2.\binom{q}{n}-1 \big).\log\big(
\binom{n+N}{n}^2.\binom{q}{n}\big)\cdot\sqrt{ p_0}}{p_0}\notag\\
&<\log(1+\frac{\epsilon}{2MN}).
\end{align*}
This is a contradiction.

We now fix a positive integer $p$ satisfying (\ref{new3}).

\noindent By (\ref{eq39}), we have
\begin{align}
 &\int\limits_{S(r)}\min_J \text{log} \prod_{j=1}^M |\psi_j^J(f)| \sigma
\geq \frac{1}{t_p}\int\limits_{S(r)}\min_J \text{log}
\Big(\prod_{\begin{matrix} \scriptstyle{1\leq j \leq M}\cr
\noalign{\vskip-0.10cm} \scriptstyle{1 \leq \ell \leq
t_p}\end{matrix}} |b_\ell \widetilde \psi_j^J(f)|\cdot
\prod_{\begin{matrix} \scriptstyle{t_p+1\leq u \leq t_{p+1}}\cr
\noalign{\vskip-0.15cm}
\scriptstyle{1 \leq v \leq M}\end{matrix}} |h_{uv}^J|\Big)\sigma\notag\\
&\quad - \frac{\varepsilon}{2}T_f(r) - o(T_f(r))\notag\\
&\geq \frac{1}{t_p}\int\limits_{S(r)} \min_J \text{log}
\frac{\prod\limits_{\begin{matrix} \scriptstyle{1 \leq j \leq M}\cr
\noalign{\vskip-0.1cm} \scriptstyle{1 \leq \ell \leq
t_p}\end{matrix}} |b_\ell \widetilde \psi_j^J(f)| \cdot
\prod\limits_{\begin{matrix} \scriptstyle{t_p+1 \leq u \leq
t_{p+1}}\cr \noalign{\vskip-0.1cm} \scriptstyle{1 \leq v \leq
M}\end{matrix}} |h_{uv}^J |}
{|W^\alpha_J|} \sigma \notag\\
&\quad + \frac{1}{t_p}\int\limits_{S(r)} \min_J
\text{log}|W^\alpha_J|\sigma
- \frac{\varepsilon}{2}T_f(r) - o(T_f(r))\notag
\end{align}
\begin{align}
&\geq - \frac{1}{t_p}\int\limits_{S(r)} \max_J \text{log}^+
\frac{|W^\alpha_J|} {\prod\limits_{\begin{matrix} \scriptstyle{1
\leq j \leq M}\cr \noalign{\vskip-0.1cm} \scriptstyle{1 \leq \ell
\leq t_p}\end{matrix}} |b_\ell \widetilde \psi_j^J(f)|\cdot
\prod\limits_{\begin{matrix} \scriptstyle{t_p+1 \leq u \leq
t_{p+1}}\cr \noalign{\vskip-0.1cm}
\scriptstyle{1 \leq v \leq M}\end{matrix}}|h_{uv}^J|} \sigma \notag\\
&\quad + \frac{1}{t_p}\int\limits_{S(r)} \text{log}|W^\alpha| \sigma
- \frac{\varepsilon}{2} T_f(r) - o(T_f(r))\,. \label{eq310}
\end{align}

For each $J \subset \{1,\dots,q\}$, $\# J=n$,
take $\beta_J$ a meromorphic function on $\C^m$ such that
\begin{align}
\big(\dots ; \beta_J b_\ell \widetilde \psi_j^J(f) : \dots :
\beta_J h_{uv}^J : \dots \big)   \label{70}
\end{align}
is a reduced representation of the meromorphic map
\begin{align*}
g_J := \Big(\dots  b_\ell \widetilde \psi_j^J(f) : \dots :
h_{uv}^J\dots \Big)^{\begin{matrix} \scriptstyle{t_p+1 \leq u \leq
t_{p+1}}\cr \noalign{\vskip-0.15cm} \scriptstyle{1 \leq v \leq
M}\end{matrix}}_{\begin{matrix} \scriptstyle{1\leq j \leq M}\cr
\noalign{\vskip-0.15cm} \scriptstyle{1 \leq \ell \leq
t_p}\end{matrix}}\ : \ \C^m \to \C P^{Mt_{p+1} -1}.
\end{align*}
By Lemma \ref{lem6} and Lemma \ref{lem7} we have
\begin{align*}
&\Big\Vert \int\limits_{S(r)} \text{log}^+ \frac{|W_J^\alpha|}
{\prod\limits_{\begin{matrix} \scriptstyle{1 \leq j \leq M}\cr
\noalign{\vskip-0.1cm} \scriptstyle{1 \leq \ell \leq
t_p}\end{matrix}} |b_{\ell} \widetilde \psi_j^J (f)| \cdot
\prod\limits_{\begin{matrix} \scriptstyle{t_p+1\leq u \leq
t_{p+1}}\cr \noalign{\vskip-0.15cm}
\scriptstyle{1 \leq v \leq M}\end{matrix}}|h_{uv}^J|}\sigma \\
&= \int\limits_{S(r)}\text{log}^+ \frac{\big| W^\alpha
(\dots,\beta_J b_{\ell} \widetilde \psi_j^J(f),\dots, \beta_J
h_{uv}^J,\dots)\big|} {\prod\limits_{\begin{matrix} \scriptstyle{1
\leq j \leq M}\cr \noalign{\vskip-0.1cm} \scriptstyle{1 \leq \ell
\leq t_p}\end{matrix}} |\beta_J b_{\ell} \widetilde \psi_j^J(f)|
\cdot \prod\limits_{\begin{matrix} \scriptstyle{t_p+1\leq u \leq
t_{p+1}}\cr \noalign{\vskip-0.1cm}
\scriptstyle{1 \leq v \leq M}\end{matrix}} |\beta_J h_{uv}^J|} \sigma\\
&= o(T_{g_J}(r))\,.
\end{align*}
On the other hand, by Corollary 5.7 in \cite{b7} or by
Theorem 5.2.29 of \cite{b20}, we have
\begin{align*}
T_{g_J}(r) &\leq \sum_{\begin{matrix} \scriptstyle{1 \leq \ell \leq
t_p}\cr \noalign{\vskip-0.15cm} \scriptstyle{1 \leq j \leq
M}\end{matrix}} T_{\frac{b_\ell \widetilde \psi_j^J(f)} {b_1
\widetilde \psi_1^J(f)}}(r) +\sum_{\begin{matrix}
\scriptstyle{t_p+1\leq u \leq t_{p+1}}\cr \noalign{\vskip-0.15cm}
\scriptstyle{1 \leq v \leq M}\end{matrix}}
T_{\frac{h_{uv}^J}{b_1 \widetilde \psi_1^J(f)}} (r) +O(1)\\
&\leq O\Big(\sum_{0 \leq i \leq n}T_{\frac{f_i}{f_0}}(r)\Big)
+ o(T_f(r)) = O (T_f(r))\,.
\end{align*}
Hence, for any $J \subset \{1,\dots,q\}$, $\# J = n$, we have
\begin{align*}
\Big\Vert \int\limits_{S(r)} \text{log}^+ \frac{|W_J^\alpha|}
{\prod\limits_{\begin{matrix} \scriptstyle{1 \leq j \leq M}\cr
\noalign{\vskip-0.15cm} \scriptstyle{1 \leq \ell \leq
t_p}\end{matrix}} |b_\ell \widetilde \psi_j^J(f)| \cdot
\prod\limits_{\begin{matrix} \scriptstyle{t_p+1 \leq u \leq
t_{p+1}}\cr \noalign{\vskip-0.15cm} \scriptstyle{1 \leq v \leq
M}\end{matrix}}|h_{uv}^J|}\sigma = o(T_f(r))\,.
\end{align*}
This implies that
\begin{align}
 \Big\Vert \int\limits_{S(r)}\max_J \text{log}^+
\frac{|W_j^\alpha|}{\prod\limits_{\begin{matrix} \scriptstyle{1 \leq
j \leq M}\cr \noalign{\vskip-0.25cm} \scriptstyle{1 \leq \ell \leq
t_p}\end{matrix}}|b_\ell
 \widetilde \psi_j^J(f)|\cdot
\prod\limits_{\begin{matrix} \scriptstyle{t_p+1 \leq u \leq
t_{p+1}}\cr \noalign{\vskip-0.15cm} \scriptstyle{1 \leq v \leq
M}\end{matrix}} |h_{uv}^J|}\sigma = o(T_f(r))\,. \label{eq312}
\end{align}
By (\ref{eq310}) and (\ref{eq312}) we get
 \begin{align*}
\Big\Vert \int\limits_{S(r)} \min_J \text{log}\prod_{j=1}^M
|\psi_j^J(f)|\sigma \geq\frac{1}{t_p}\int\limits_{S(r)}
\text{log}|W^\alpha| \sigma - \frac{\varepsilon}{2}T_f(r) -
o(T_f(r)).
\end{align*}
Therefore, by (\ref{eq37}), we obtain that
\begin{align}
\Big\Vert \frac{1}{A\cdot t_{p}}\int\limits_{S(r)} \text{log}
\frac{(\prod_{j=1}^q |Q_j(f)| )^{A\cdot t_{p}}}{|W^\alpha| }\sigma
&\geq (q-n-1)d\cdot T_f(r) \notag\\
 -(\frac{\epsilon}{2}+\frac{\epsilon}{2A})\cdot T_f(r) - o(T_f(r)).
\label{13}
\end{align}

We recall that
 \begin{align*}
Q_j = \sum_{I \in \Cal T_d} a_{jI}x^I,\quad a_{jI} \in \Cal K_f \quad
(j = 0,\dots,q).
\end{align*}
Let again $T = (\dots,t_{kI},\dots)$ \ ($k \in \{0,\dots,q\}$, $I \in \Cal T_d$)
be a family of variables and
\begin{align*}
\widetilde Q_j = \sum_{I \in \Cal T_d} t_{jI}x^I \in \Z[T,x],\quad
j = 0,\dots, q.
\end{align*}
For each $H\subset \{1,...,q\},\#H=(n+1),$ let $\widetilde R_H \in \Z [T]$ be the resultant of $\{\widetilde Q_j\}_{j \in  H}.$ By Proposition 2.1 there exists a positive integer $s$ (without loss of generality) common for all $H$ and polynomials $\{ \widetilde b_{ij}^H\} (0\leq i\leq n; j\in H)$ in $\Z[T,x],$ which are zero or  homogeneous in $x$ of degree $(s-d),$ such that
\begin{align*}
x_i^s \cdot \widetilde  R_H = \sum_{j\in H} \widetilde b_{ij}^H\cdot \widetilde  Q_j\quad
\text{for all \ }  i \in \{0,\dots,n\},
\end{align*}
and $R_H = \widetilde R_H(\dots, a_{kI}, \dots) \not\equiv 0$.  We note again that $ R_H\in  \Cal K_f.$

\noindent Set
\begin{align*}
b_{ij}^H := \widetilde b_{ij}^H\big((\dots,a_{jI},\dots), (x_0,\dots,x_n)\big).\quad
\end{align*}
Then we have
\begin{align*}
f_i^s \cdot  R_H = \sum_{j\in H}b_{ij}^H (f)\cdot  Q_j(f)\quad
\text{for all}\  i \in \{0,\dots,n\}.
\end{align*}
This implies, since $(f_0:\cdots :f_n)$ is a reduced representation, i.e. $\{f_0=\cdots=f_n=0\}$  of codimension $\geq 2$,  that
\begin{align}
\nu_{R_H}\geq \min\limits_{j\in
H}\nu_{Q_j(f)}-\sum\limits_{\begin{matrix} \scriptstyle{0 \leq i
\leq n}\cr \noalign{\vskip-0.15cm} \scriptstyle{j\in H}\end{matrix}}
\nu_{\frac{1}{b_{ij}^H(f)}}
\end{align}
outside an analytic subset of codimension $\geq 2.$
Then, we have
\begin{align}
\nu_{\frac{(\prod_{j=1}^q Q_j(f))^{A\cdot t_{p}}}{W^\alpha }}&\leq\max_{J\subset  \{1, ... , q\}, \#J=n}\nu_{\frac{(\prod_{j\in J}Q_j(f))^{ A\cdot t_{p}}}{W^\alpha }}\notag\\
&\quad+ \min_{E  \subset  \{1, ... , q\}, \# E=q-n}\nu_{(\prod_{j \in  E}Q_j(f))^{A\cdot t_{p}}}\notag\\
&\leq\max_{J\subset  \{1, ... , q\}, \#J=n}\nu_{\frac{(\prod_{j\in J}Q_j(f))^{A\cdot t_{p}}}{W_J^\alpha} }\notag\\
&\quad+ A t_{p}\min_{E  \subset  \{1, ... , q\}, \# E=q-n}\nu_{\prod_{j\in E}Q_j(f)}\notag\\
&\leq\max_{J\subset  \{1, ... , q\}, \#J=n}\nu_{\frac{(\prod_{j\in  J}Q_j(f))^{A\cdot t_{p}}}{W_J^\alpha }}\notag\\
&\quad+(q-n) A t_{p}\sum\limits_{\begin{matrix} \scriptstyle{H \in
\{1, ...,q\}}\cr \noalign{\vskip-0.15cm} \scriptstyle{\#
H=n+1}\end{matrix}} \Big (\nu_{R_H}+ \sum\limits_{\begin{matrix}
\scriptstyle{0 \leq i \leq n}\cr \noalign{\vskip-0.15cm}
\scriptstyle{j\in H}\end{matrix}} \nu_{\frac{1}{b_{ij}^H(f)}} \Big )
\label{61}
\end{align}
outside an analytic subset of codimension $\geq 2$
(note that
\begin{align*}
W^\alpha_{J} = W^\alpha(b_\ell \widetilde \psi_j^{J}(f), 1 \leq \ell
\leq t_{p},
1 \leq j \leq M, h_{uv}^J, t_{p}+1 \leq u \leq t_{p+1}, 1 \leq v \leq M)\\
= W^\alpha(b_\ell (c^{J}_{j I^{J}_j})^{-1} \psi_j^{J}(f), 1 \leq
\ell \leq t_{p},
1 \leq j \leq M, h_{uv}^J, t_{p}+1 \leq u \leq t_{p+1}, 1 \leq v \leq M)\\
=\text{det}C_{J}\cdot W^\alpha, \quad C_J\in GL(t_{p+1}M, \C)).
\end{align*}
Take $L=Mt_{p+1}-1\leq \binom{n+N}{n}t_{p_0+1}-1$ (note that $p\leq
p_0).$ We recall that
\begin{align*}
N&=d\cdot[2(n+1)(2^n-1)(nd+1)\epsilon^{-1}+n+1],\\
p_0&=\big[\frac{\big( \binom{n+N}{n}^2.\binom{q}{n}-1
\big).\log\big(
\binom{n+N}{n}^2.\binom{q}{n}\big)}{\log(1+\frac{\epsilon}{2\binom{n+N}{n}N})}+1\big]^2,\;\text{and}\\
 t_{p_0+1}
&\overset{(\ref{new4})}{<}\Bigg(\binom{n+N}{n}^2.\binom{q}{n}+p_0\Bigg)^{
\binom{n+N}{n}^2.\binom{q}{n}-1}.
\end{align*}
 Furthermore, in the case of fixed
hypersurfaces $(Q_j\in\C[x_0,\dots,x_n], j=1,\dots, q),$  that is in the case of  
$(Q_j\in\C[x_0,\dots,x_n], j=1,\dots, q)$ with constant coefficients $a_{jI} \in \C$,
 by the
definition of $\mathcal L(p)$
we have $t_p=1$ for all positive
integer $p,$ and then $L=M-1.$

 For each $J \subset
\{1,\dots,q\}$, $\# J=n$, we use again the reduced representation
(see (\ref{70}))
\begin{align*}
\big(\dots ; \beta_J b_\ell \widetilde \psi_j^J(f) : \dots :
\beta_J h_{uv}^J : \dots \big)
\end{align*}
 of the meromorphic map
\begin{align*}
g_J := \Big(\dots  b_\ell \widetilde \psi_j^J(f) : \dots :
h_{uv}^J\dots \Big)^{\begin{matrix} \scriptstyle{t_p+1 \leq u \leq
t_{p+1}}\cr \noalign{\vskip-0.15cm} \scriptstyle{1 \leq v \leq
M}\end{matrix}}_{\begin{matrix} \scriptstyle{1\leq j \leq M}\cr
\noalign{\vskip-0.15cm} \scriptstyle{1 \leq \ell \leq
t_p}\end{matrix}}\ : \ \C^m \to \C P^{Mt_{p+1} -1}.
\end{align*}
For any $J=(j_1,\dots, j_n)\subset\{1,\dots,q\},$ set
\begin{align*}
\nu_J:=\nu_{\beta_J}+\nu_{\frac{1}{\beta_J}}+\sum\limits_{\begin{matrix}
\scriptstyle{1 \leq j \leq M}\cr \noalign{\vskip-0.15cm}
\scriptstyle{1 \leq \ell \leq t_p}\end{matrix}} \nu_{\frac{1}{b_\ell
(c^{J}_{j I^{J}_j})^{-1}\gamma_j^J(f) }}+\sum_{\begin{matrix}
\scriptstyle{t_p+1 \leq u \leq t_{p+1}}\cr \noalign{\vskip-0.15cm}
\scriptstyle{1 \leq v \leq M}\end{matrix}}\nu_{\frac{1}{h_{uv}^{J}}}
\end{align*}
(note that $\gamma_j^J$ is defined by (\ref{new10})). It is easy to
see that $ N_{\nu_J}(r)=o(T_f(r)).$

By Proposition \ref{lemy} and Lemmas \ref{lemz}-\ref{lem6}, we have
(outside an analytic subset of codimension $\geq 2)$
\begin{align*}
&\nu_{\frac{\big(\prod\limits_{j\in J}Q_j(f)\big)^{A\cdot
t_{p}}}{W^\alpha_{J}}}+\sum\limits_{\begin{matrix} \scriptstyle{1
\leq j \leq M}\cr \noalign{\vskip-0.15cm} \scriptstyle{1 \leq \ell
\leq t_p}\end{matrix}}\nu_{b_\ell (c^{J}_{j I^{J}_j})^{-1}
\gamma_j^J(f) } +\\
&+ \sum\limits_{\begin{matrix} \scriptstyle{t_p+1 \leq
u \leq t_{p+1}}\cr \noalign{\vskip-0.15cm} \scriptstyle{1 \leq v
\leq M}\end{matrix}}\nu_{
h_{uv}^{J}}-\nu_{\big(\frac{(\prod\limits_{j\in J}Q_j(f))^{A\cdot
t_{p}}}{W^\alpha_{J}}\big)^{-1}}-\nu_J\notag\\
&\leq \nu_{\frac{\big(\prod\limits_{\begin{matrix} \scriptstyle{1 \leq j
\leq M}\cr \noalign{\vskip-0.15cm} \scriptstyle{1 \leq \ell \leq
t_p}\end{matrix}}b_\ell (c^{J}_{j
I^{J}_j})^{-1}\big)\cdot\big(\prod\limits_{j\in
J}Q_j(f)\big)^{A\cdot t_{p}}
\cdot\prod\limits_{\ell=1}^M(\gamma_\ell^J(f))^{t_p}\cdot\prod\limits_{\begin{matrix}
\scriptstyle{t_p+1 \leq u \leq t_{p+1}}\cr \noalign{\vskip-0.15cm}
\scriptstyle{1 \leq v \leq M}\end{matrix}}h_{uv}^{J}}{W^\alpha_{J}}}\\
&\overset{(\ref{45})}{=}\nu_{\frac{\prod\limits_{\begin{matrix}
\scriptstyle{1 \leq j \leq M}\cr \noalign{\vskip-0.15cm}
\scriptstyle{1 \leq \ell \leq t_p}\end{matrix}} \big(b_\ell
(c^{J}_{j I^{J}_j})^{-1} \psi_j^{J}(f)\big)
\prod\limits_{\begin{matrix} \scriptstyle{t_p+1 \leq u \leq
t_{p+1}}\cr \noalign{\vskip-0.15cm} \scriptstyle{1 \leq v \leq
M}\end{matrix}}h_{uv}^{J}}{W^\alpha_{J}}}
\\ & 
=
\nu_{\frac{\prod\limits_{\begin{matrix} \scriptstyle{1 \leq j \leq
M}\cr \noalign{\vskip-0.15cm} \scriptstyle{1 \leq \ell \leq
t_p}\end{matrix}}\big(\beta_J b_\ell (c^{J}_{j I^{J}_j})^{-1}
\psi_j^{J}(f)\big)  \prod\limits_{\begin{matrix} \scriptstyle{t_p+1
\leq u \leq t_{p+1}}\cr \noalign{\vskip-0.15cm}
\scriptstyle{1 \leq v \leq M}\end{matrix}}\beta_J h_{uv}^{J}}{W^\alpha (\dots, \beta_J b_\ell (c^{J}_{j I^{J}_j})^{-1} \psi_j^{J}(f), \dots , \beta_J h_{uv}^{J}, \dots )}}\notag\\
&\leq \sum\limits_{\begin{matrix} \scriptstyle{1 \leq j \leq M}\cr
\noalign{\vskip-0.15cm} \scriptstyle{1 \leq \ell \leq
t_p}\end{matrix}} \text{min}\{\nu_{\beta_Jb_\ell (c^{J}_{j
I^{J}_j})^{-1} \psi_j^{J}(f) },L\} +\sum\limits_{\begin{matrix}
\scriptstyle{t_p+1 \leq u \leq t_{p+1}}\cr \noalign{\vskip-0.15cm}
\scriptstyle{1 \leq
v \leq M}\end{matrix}}\text{min}\{\nu_{\beta_J h_{uv}^{J}},L\}
%\notag\\
\end{align*}
\begin{align*}
&\overset{(\ref{new10})}{\leq}
At_p\sum_{i=1}^n\min\{\nu_{Q_{j_i}(f)},L\} +
\sum\limits_{\begin{matrix} \scriptstyle{1 \leq j \leq M}\cr
\noalign{\vskip-0.15cm} \scriptstyle{1 \leq \ell \leq
t_p}\end{matrix}} \nu_{\beta_Jb_\ell (c^{J}_{j I^{J}_j})^{-1}
\gamma_j^J(f) }+\sum\limits_{\begin{matrix} \scriptstyle{t_p+1 \leq
u \leq t_{p+1}}\cr \noalign{\vskip-0.15cm} \scriptstyle{1 \leq v
\leq
M}\end{matrix}}\nu_{\beta_J h_{uv}^{J}}
\notag\\
&\leq At_p\sum_{j=1}^q\min\{\nu_{Q_{j}(f)},L\} +
\sum\limits_{\begin{matrix} \scriptstyle{1 \leq j \leq M}\cr
\noalign{\vskip-0.15cm} \scriptstyle{1 \leq \ell \leq
t_p}\end{matrix}}\nu_{b_\ell (c^{J}_{j I^{J}_j})^{-1} \gamma_j^J(f)
} +\sum\limits_{\begin{matrix} \scriptstyle{t_p+1 \leq u \leq
t_{p+1}}\cr \noalign{\vskip-0.15cm} \scriptstyle{1 \leq v \leq
M}\end{matrix}}\nu_{ h_{uv}^{J}}+O(\nu_{\beta_J}).
\end{align*}
 This implies that for any  J (outside an analytic subset of codimension $\geq
 2$)
\begin{align}
\nu_{\frac{(\prod_{j\in J} Q_j(f))^{A\cdot t_{p}}}{W_J^\alpha
}}&\leq At_p\sum_{j=1}^q\min\{\nu_{Q_{j}(f)},L\}
+\nu_{\big(\frac{(\prod\limits_{j\in J}Q_j(f))^{A\cdot
t_{p}}}{W^\alpha_{J}}\big)^{-1}}+\nu\,,\label{nx}
\end{align}
where $\nu$ is a divisor on $\C^m$ such that $N_\nu(r)=o(T_f(r)).$
\vspace{0.5cm}

\noindent Since the set $\{ \overline{\nu_{\big(\frac{(\prod\limits_{j\in J}Q_j(f))^{A\cdot
t_{p}}}{W^\alpha_{J}}\big)} >0 \}} \cap 
\overline{\{ \nu_{\big(\frac{(\prod\limits_{j\in J}Q_j(f))^{A\cdot
t_{p}}}{W^\alpha_{J}}\big)^{-1}}>0\}}$ is an analytic set of codimension $\geq 2$, and the terms on the right hand side
of (\ref{nx}) are all nonnegative, this implies that for all J (outside an analytic subset of codimension $\geq 2$)

\begin{align*}
\nu_{\frac{(\prod_{j\in J} Q_j(f))^{A\cdot t_{p}}}{W_J^\alpha
}}&\leq At_p\sum_{j=1}^q\min\{\nu_{Q_{j}(f)},L\}
+\nu\,.
\end{align*}

\noindent Then, by (\ref{61}) we have
\begin{align}
\nu_{\frac{(\prod_{j=1}^q Q_j(f))^{A\cdot t_{p} }}{W^\alpha }}&\leq
At_p\sum_{j=1}^q\min\{\nu_{Q_{j}(f)}, L\}
+\nu \notag\\
&\quad\quad\quad+(q-n) A t_{p}\sum\limits_{\begin{matrix}
\scriptstyle{H \subset \{1, ...,q\}}\cr \noalign{\vskip-0.15cm}
\scriptstyle{\# H=n+1}\end{matrix}}\Big (
\nu_{R_H}+\sum\limits_{\begin{matrix} \scriptstyle{0 \leq i \leq
n}\cr \noalign{\vskip-0.15cm} \scriptstyle{j\in H}\end{matrix}}
\nu_{\frac{1}{b_{ij}^H(f)}} \Big ) \label{64}
\end{align}
outside an analytic subset of codimension $\geq 2.$

 By Jensen's formula and by (\ref{64}), we get
\begin{align*}
&\int \limits_{S(r)} \text{log} \frac{(\prod_{j=1}^q |Q_j(f)|
)^{A\cdot t_{p}}}{|W^\alpha| }\sigma\\
& =N_{\frac{(\prod_{j=1}^q Q_j(f)
)^{A\cdot t_{p}}}{W^\alpha}}(r) -N_{\Big(\frac{(\prod_{j=1}^q Q_j(f)
)^{A\cdot t_{p}}}{W^\alpha}\Big)^{-1}}(r) +O(1)
\\
&\leq A t_{p}\sum_{j=1}^q N^{(L)}_f(r,Q_j)+N_{\nu}(r) +O(1)\\
&\quad\quad\quad\qquad+(q-n) A t_{p}\sum\limits_{\begin{matrix}
\scriptstyle{H \subset \{1, ...,q\}}\cr \noalign{\vskip-0.15cm}
\scriptstyle{\# H=n+1}\end{matrix}}\Big (
N_{R_H}(r)+\sum\limits_{\begin{matrix} \scriptstyle{0 \leq i \leq
n}\cr \noalign{\vskip-0.15cm}
\scriptstyle{j\in H}\end{matrix}} N_{\frac{1}{b_{ij}^H(f)}}(r) \Big )\\
&\leq A t_{p}\sum_{j=1}^q  N^{(L)}_f(r,Q_j)+o(T_f(r))
\end{align*}
(note that $R_H\in\mathcal K_f,\; b_{ij}^H\in\mathcal
K_f[x_0,\dots,x_n]).$

\noindent Combining with (\ref{13}), we have
\begin{align}
\Vert (q-n-1-\varepsilon) T_f(r) \leq \sum_{j=1}^q \frac{1}{d}
N^{(L)}_f(r,Q_j),\label{new8}
\end{align}
(note that $A>1$).

We now prove the theorem for the general case: $\deg Q_j=d_j.$
Denote by  $d$ the least common multiple of $d_1,\dots, d_q$ and put
$d_j^*:=\frac{d}{d_j}.$ By (\ref{new8}) with the moving
hypersurfaces $Q_j^{d_j^*}$ $(j\in\{1\dots,q\})$ of common degree
$d,$ we have
\begin{align*}
\Vert (q-n-1-\varepsilon) T_f(r) \leq \sum_{j=1}^q \frac{1}{d}
N^{(L)}_f(r,Q_j^{d_j^*})\leq\sum_{j=1}^q \frac{d_j^*}{d}
N^{([\frac{L}{d_j^*}+1])}_f(r,Q_j)\\
\leq \sum_{j=1}^q \frac{1}{d_j} N^{(L_j)}_f(r,Q_j),
\end{align*}
where $L_j:=[\frac{d_jL}{d}+1].$ This completes the proof of  the Main
Theorem. $\square$

\vspace{1cm}

 \noindent  Gerd Dethloff$^{1-2} $ \\
 $^1$ Universit\'e Europ\'eenne de Bretagne, France\\
 $^2$
Universit\'{e} de Brest \\
   Laboratoire de math\'{e}matiques \\
UMR CNRS 6205\\
6, avenue Le Gorgeu, BP 452 \\
   29275 Brest Cedex, France \\
e-mail: gerd.dethloff@univ-brest.fr\\

\noindent Tran Van Tan\\
Department of Mathematics\\
  Hanoi National University of Education\\
 136-Xuan Thuy street, Cau Giay, Hanoi, Vietnam\\
e-mail: tranvantanhn@yahoo.com

\end{document}